\documentclass[english]{smfart}
\usepackage[all]{xypic}
\usepackage{amssymb}
\usepackage{amsfonts}
\usepackage{amscd}
\usepackage[all]{xypic}
\usepackage[english]{babel}

\swapnumbers

\newtheorem{theorem}[subsection]{Theorem}
\newtheorem{lem}[subsection]{Lemma}
\newtheorem{cor}[subsection]{Corollary}
\newtheorem{prop}[subsection]{Proposition}

\theoremstyle{definition}
\newtheorem{definition}[subsection]{Definition}

\newtheorem{def-prop}[subsection]{Definition-Proposition}

\theoremstyle{remark}
\newtheorem{remark}[subsection]{Remark}

\theoremstyle{plain}

\def\boxit#1#2{\setbox1=\hbox{\kern#1{#2}\kern#1}%
\dimen1=\ht1 \advance\dimen1 by #1
\dimen2=\dp1 \advance\dimen2 by #1
\setbox1=\hbox{\vrule height\dimen1 depth\dimen2\box1\vrule}%
\setbox1=\vbox{\hrule\box1\hrule}%
\advance\dimen1 by .4pt \ht1=\dimen1
\advance\dimen2 by .4pt \dp1=\dimen2 \box1\relax}

\def\ie{{\emph {i.e. \/}}}

\def\AA{{\mathbf A}}
\def\BB{{\mathbf B}}
\def\CC{{\mathbf C}}

\def\LL{{\mathbf L}}

\def\NN{{\mathbf N}}

\def\PP{{\mathbf P}}
\def\QQ{{\mathbf Q}}
\def\RR{{\mathbf R}}

\def\ZZ{{\mathbf Z}}

\def\cJ{{\mathcal J}}

\def\cL{{\mathcal L}}
\def\cM{{\mathcal M}}

\def\cO{{\mathcal O}}

\def\cV{{\mathcal V}}

\mathchardef\alphag="7C0B
\mathchardef\betag="7C0C
\mathchardef\gammag="7C0D
\mathchardef\deltag="7C0E
\mathchardef\varepsilong="7C22
\mathchardef\varphig="7C27
\mathchardef\psig="7C20
\mathchardef\zetag="7C10
\mathchardef\epsilong="7C0F
\mathchardef\rhog="7C1A
\mathchardef\taug="7C1C
\mathchardef\upsilong="7C1D
\mathchardef\iotag="7C13
\mathchardef\thetag="7C12
\mathchardef\pig="7C19
\mathchardef\sigmag="7C1B
\mathchardef\etag="7C11
\mathchardef\omegag="7C21
\mathchardef\kappag="7C14
\mathchardef\lambdag="7C15
\mathchardef\mug="7C16
\mathchardef\xig="7C18
\mathchardef\chig="7C1F
\mathchardef\nug="7C17
\mathchardef\varthetag="7C23
\mathchardef\varpig="7C24
\mathchardef\varrhog="7C25
\mathchardef\varsigmag="7C26
\mathchardef\Omegag="7C0A
\mathchardef\Thetag="7C02
\mathchardef\Sigmag="7C06
\mathchardef\Deltag="7C01
\mathchardef\Phig="7C08
\mathchardef\Gammag="7C00
\mathchardef\Psig="7C09
\mathchardef\Lambdag="7C03
\mathchardef\Xig="7C04
\mathchardef\Pig="7C05
\mathchardef\Upsilong="7C07

\DeclareMathOperator*{\Spec}{Spec}
\def\ord{{\rm ord}}

\begin{document}

\title[Motivic integration and quotient singularities]{Motivic
integration,
quotient singularities
and the McKay correspondence}

\author{Jan Denef}
\address{University of Leuven, Department of Mathematics,
Celestijnenlaan 200B, 3001 Leu\-ven, Bel\-gium }
\email{ Jan.Denef@wis.kuleuven.ac.be}
\urladdr{http://www.wis.kuleuven.ac.be/wis/algebra/denef.html}

\author{Fran\c cois Loeser}

\address{Centre de Math{\'e}matiques,
Ecole Polytechnique,
F-91128 Palaiseau
(UMR 7640 du CNRS), {\rm and}
Institut de Math{\'e}matiques,
Universit{\'e} P. et M. Curie, Case 82,
4 place Jussieu,
F-75252 Paris Cedex 05
(UMR 7596 du CNRS)}
\email{loeser@math.polytechnique.fr}
\urladdr{http://math.polytechnique.fr/cmat/loeser}





\maketitle

\section*{Introduction}
Let $X$ be an algebraic variety, not necessarily smooth,
over a field $k$ of
characteristic
zero. We denote by $\cL (X)$ the $k$-scheme of formal arcs on $X$~: 
$K$-points of $\cL (X)$ correspond
to formal arcs $ {\rm Spec} \, K [[t]] \rightarrow X$,
for $K$ any field containing $k$.
In a recent paper \cite{Arcs}, we developped an integration theory on
the space $\cL (X)$
with values in $\widehat \cM$, a certain ring
completion 
of the Grothendieck ring $\cM$
of algebraic varieties over $k$ (the definition
of these rings is recalled in section
\ref{GK}), based on
ideas of Kontsevich \cite{K}. In the most interesting cases, the integrals we
consider
belong to a much smaller ring
$\overline \cM_{\rm loc} [(\frac{\LL - 1}{\LL^{i} - 1})_{i \geq 1}]$,
on which the usual Euler characteristic and
Hodge polynomial may be extended in a natural way
to an Euler characteristic ${\rm Eu}$ and a Hodge polynomial $H$
belonging respectively to $\QQ$
and the ring 
$\ZZ [u, v][(uv)^{-1}][(\frac{uv - 1}{(uv)^{i} - 1})_{i \geq
1}]$. When $X$ is smooth and one considers
the total measure of $\cL (X)$,
these invariants reduce to the usual
Euler characteristic and 
Hodge polynomial, but in general one obtains interesting new invariants
(see \cite{BD},\cite{BB},\cite{Ba2},\cite{Arcs},\cite{Veys}).

\medskip

When $X$ is a normal variety with at most Gorenstein canonical
singularities, one can use the canonical class to define a measure
$\mu^{Gor} (A)$ for certain subsets $A$ of $\cL (X)$.
Now assume $X$ is 
the quotient of the affine space
$\AA^{n}_{k}$ by a finite subgroup $G$ of order $d$
of
${\rm SL}_{n} (k)$. We make the assumption $k$ contains all $d$-th
roots of unity.
We denote by  $\cL (X)_{0}$
the set of arcs whose origin is the point
0 in $X$.
One of the
main results of the present paper is Theorem \ref{MT},
which expresses
$\mu^{Gor} (\cL (X)_{0})$ in terms of representation theoretic
weights $w (\gamma)$
of the conjugacy classes of elements $\gamma$ of $G$,
defined as
$w (\gamma) := \sum_{1 \leq i \leq n} e_{\gamma, i} / d$,
with
$1 \leq e_{\gamma, i} \leq d$
and $\xi^{e_{\gamma, i}}$
the eigenvalues of $\gamma$ for $ i = 1, \ldots, n$,
$\xi$ being a fixed primitive $d$-th root of unity 
in $k$.
More precisely, the image of
$\mu^{Gor} (\cL (X)_{0})$ is equal to that of
$\sum_{[\gamma] \in {\rm Conj} (G)} \, \LL^{- w (\gamma)}$
in a certain quotient $\widehat \cM_{/} $ of
$\widehat \cM$,
with $\LL$ the class of the affine line. The 
quotient $\widehat \cM_{/} $ is defined by requiring that the class of
a quotient of a vector space $V$ by a finite group acting linearly should
be that of $V$. This condition is mild enough to guarantee that
$\mu^{Gor} (\cL (X)_{0})$ and 
$\sum_{[\gamma] \in {\rm Conj} (G)} \, \LL^{- w (\gamma)}$
have the same image in 
$\widehat
K_{0} ({\rm CHM}_{k})$, an appropriate completion of $K_{0} ({\rm
CHM}_{k})$,
the Grothendieck group of the pseudo-abelian
category
of Chow motives over $k$, and in particular have the
same Hodge
polynomial and Euler characteristic.
This result - at least for the  Hodge realization -
is due to Batyrev
\cite{Ba} and implies, when
$X$ has a crepant resolution,
a
form of the McKay correspondence
which has been conjectured by Reid \cite{Reid}
and proved by Batyrev \cite{Ba}.

\medskip

The aim of the present paper is to present an alternative proof of 
Batyrev's result and also to develop further
some basic properties of motivic integration which were not covered in
\cite{Arcs}.
Though Batyrev also uses integration on spaces of arcs,
the approach we follow here, which was
inspired to us by Maxim Kontsevich, is somewhat different.
One of the main differences is that we are able to
work
directly on the singular space $X$
instead of going to desingularisations. 
This allows us to have a more local approach,
in the sense that we can directly calculate
the part of the motivic integral
coming from each conjugacy class. More precisely,
for each element $\gamma$ in the group $G$, we consider
$\cL (X)^{g}_{0, \gamma}$, the set of arcs $\varphi$
in $\cL (X)_0$, which are
not
contained
in the discriminant and may be lifted in $\cL (\AA^{n}_{k})$ to a fractional
arc
$\tilde \varphi (t^{1 / d})$ such that
$\tilde \varphi (\xi t^{1 / d}) = \gamma \, \tilde \varphi (t^{1 /
  d})$.
We prove that the image of 
$\mu^{Gor} (\cL (X)^{g}_{0, \gamma})$ 
in $\widehat \cM_{/}$
is equal to that of
$\LL^{- w (\gamma)}$.

\medskip

Let us now briefly review the content of the paper.
In section \ref{pre}, we recall some material
on semi-algebraic geometry over $k ((t))$ from
\cite{Arcs}. In fact, we need to generalize slightly
semi-algebraic geometry as developped in 
\cite{Arcs} to ``$k [t]$-semi-algebraic geometry'' which allows 
expressions involving $t$, since $k [t]$-morphisms
naturally appear in section \ref{so}.
Fortunately this is quite harmless, since most proofs remain the same.
This material on $k [t]$-semi-algebraic geometry might be useful elsewhere.
One of the main technical difficulties of the section is 
Theorem \ref{CV}
were we  extend
the crucial
change of variables formula \cite{Arcs} to certain maps which are not
birational.
Section \ref{so} is the heart of the 
paper, namely the study of the local action of the group $G$ on arcs.
We are then able to deduce the main results in section
\ref{mq}.
In section \ref{mot} we explain how 
one can deduce statements at the level of Chow motives and then
realizations,
and in section \ref{mk} we express the main results in terms of 
resolutions of singularities and we explain the relation with
McKay's correspondence.

\medskip

Let us remark that $k [t]$-semi-algebraic sets appear quite naturally in the
problem, since the set
$\cL (X)^{g}_{0, \gamma}$ is $k [t]$-semi-algebraic.
Nevertheless, it is possible to avoid the use
of $k [t]$-semi-algebraic geometry here, by using properties of measurable
sets which are developped in the appendix, in particular the fact,
proved in Theorem \ref{image}, that
the image of a measurable set under a $k [t]$-morphism, for varieties of
the same dimension,  is again measurable.

\bigskip

{\it Acknowledgements}~: We thank
Maxim Kontsevich for a conversation which greatly inspired
the approach followed in the present paper.
We thank also Eduard Looijenga who helped us to improve 
the paper,   and in particular 
to avoid the use of motives in
Lemma \ref{4.5}.


\renewcommand{\theequation}{\thesubsection.\arabic{equation}}

\section{Preliminaries on semi-algebraic geometry and\\
motivic integration}\label{pre}\subsection{}In the present paper
by a variety over $k$, or variety,  we always
mean 
a reduced separated scheme of finite type over a field $k$ that will 
be assumed to be of characteristic zero throughout the paper.
If $X$ is a variety, we shall denote by
$X_{\rm sing}$ the singular locus of $X$.

\subsection{}\label{defarrcs}For $X$ a variety over $k$, we will denote
by $\cL (X)$ the {\emph {scheme of 
germs of
arcs on $X$}}. It is a scheme over $k$ and for any field extension
$k \subset K$ there is a natural bijection
$$\cL (X) (K) \simeq {\rm Mor}_{k-{\rm schemes}} (\Spec K [[t]], X)
$$ 
between the set of $K$-rational points of $\cL (X)$
and the set of 
germs of arcs with coefficients in  $K$ on $X$.
We will call
$K$-rational points of $\cL (X)$, for $K$
a field extension of $k$, arcs on $X$, and $\varphi (0)$ will be called 
the origin of the arc $\varphi$.
More precisely the
scheme $\cL (X)$ is defined as the projective limit
$$
\cL (X) := \varprojlim \cL_{n} (X)
$$
in the category of $k$-schemes of the schemes
$\cL_{n}(X)$ representing the functor
$$R \mapsto {\rm Mor}_{k-{\rm schemes}} (\Spec R [t] / t^{n+1} R[t], X)$$
defined on the category of $k$-algebras. (The existence of $\cL_{n}(X)$
is well known (cf. \cite{Arcs}) and the projective limit exists since 
the transition morphisms are affine.)
We shall denote by $\pi_{n}$ the canonical morphism, corresponding to 
truncation of arcs,
$$
\pi_{n} : \cL (X) \longrightarrow \cL_{n} (X).
$$
The schemes $\cL (X)$ and $\cL_{n} (X)$ will always be
considered with their reduced structure.
If $W$ is a subscheme of
$X$, we set $\cL (X)_{W} = \pi_{0}^{-1} (W)$.

Since, in section \ref{so}, we shall lift arcs to Galois covers,
we also have to consider ``ramified'' arcs, so we define 
similarly, for $d \geq 1$ an integer, 
the scheme $\cL^{1 / d} (X)$ as the projective limit
$$
\cL^{1 / d} (X) := \varprojlim \cL_{n}^{1 / d} (X)
$$
in the category of $k$-schemes of the schemes
$\cL_{n}^{1 / d}(X)$ representing the functor
$$R \mapsto {\rm Mor}_{k-{\rm schemes}} (\Spec R [t^{1 / d}] /
t^{(n+1)
/ d} R[t^{1 / d}], X)$$
defined on the category of $k$-algebras. Of course
the schemes $\cL^{1 / d} (X)$ are all isomorphic to
$\cL (X)$.
We shall still denote by $\pi_{n}$ the canonical morphism
$$
\pi_{n} : \cL (X)^{1 / d} \longrightarrow \cL_{n}^{1 / d} (X).
$$
and for $W$ a subscheme of
$X$, we set $\cL^{1 / d} (X)_{W} = \pi_{0}^{-1} (W)$.

The above definitions extend to the case where
$X$ is a reduced and separated scheme of finite type over $k
[t]$. For $n$ in $\NN$, one defines the 
$k$-scheme
$\cL_{n} (X)$ as  representing
the functor
$$
R \mapsto {\rm Mor}_{k[t]-{\rm schemes}} ({\rm Spec}
\, R [t]/t^{n+1} R [t], X),
$$
defined on the category of $k$-algebras, and one sets
$\cL (X) := \varprojlim \cL_{n} (X)$.
The existence of $\cL_{n} (X)$ is well known, cf. \cite{B-L-R} p.276,
and again
the
projective limit exists
since the transition morphisms are affine. We shall still denote 
by $\pi_n$ the canonical morphism $\cL (X) \rightarrow \cL_n (X)$.

\subsection{}Let $X$ and $Y$ be $k$-varieties.
A function 
$h : \cL (Y) \rightarrow \cL (X )$
will be called a {\emph {$k [t]$-morphism}}
if it is
induced by a morphism of $k[t]$-schemes
$Y \otimes_{k} k [t] \rightarrow X \otimes_{k} k [t]$. 
We shall denote by the same symbol a $k [t]$-morphism
and the corresponding morphism of $k[t]$-schemes.

\subsection{}We now introduce the concept of
semi-algebraic and 
$k [t]$-semi-algebraic subsets of the space of arcs $\cL (X)$. 
The main motivation for introducing such objects is that in general 
being
a subset of $\cL (X)$ defined by (boolean combination of)
algebraic conditions is not a property which is conserved by taking
images, \ie 
Theorem \ref{Pas}
and 
Proposition \ref{Pas2} (1) would not remain true when replacing 
``semi-algebraic'' by ``(boolean combination of)
algebraic''.

From now on we will denote by
$\bar k$ a fixed algebraic closure of $k$, and by $\bar k((t))$
the fraction field of $\bar k [[t]]$, where $t$ is one variable.
Let $x_1, \ldots, x_m$
be variables running over $\bar k ((t))$ and let $\ell_1, \ldots,
\ell_r$ be variables running over $\ZZ$.
A {\it semi-algebraic} (resp.
$k [t]$-{\it semi-algebraic})
condition
$\theta (x_1, \ldots, x_m; \ell_1, \ldots, \ell_r)$
is a finite boolean combination of conditions of the form
\begin{align*}
\text{(1)}&&\ord_t f_1 (x_1, \ldots, x_m) \geq
\ord_t f_2 (x_1, \ldots, x_m) + L (\ell_1, \ldots, \ell_r)\\
\text{(2)}&&\ord_t f_1 (x_1, \ldots, x_m) \equiv
L (\ell_1, \ldots, \ell_r) \pmod d\\
\intertext{and}
\text{(3)}&&h (\overline{ac} (f_1 (x_1, \ldots, x_m)),
\ldots,
\overline{ac} (f_{m'} (x_1, \ldots, x_m))) = 0,
\end{align*}
where $f_i$ are polynomials with coefficients in  $k$
(resp. $f_i$ are
polynomials with coefficients in  $k [t]$), 
$h$ is a polynomial with coefficients in  $k$,
$L$ is a polynomial of degree $\leq 1$ over $\ZZ$, $d \in \NN$,
and $\overline{ac} (x)$ is the coefficient of lowest degree
of $x$ in $\bar k ((t))$ if $x \not= 0$, and is equal to 0 otherwise.
Here we use the convention that $\infty + \ell =  \infty$
and
$\infty \equiv \ell \; {\rm mod} \, d$, for all $\ell \in \ZZ$.
In particular, the algebraic (resp. $k [t]$-algebraic)
condition 
$f (x_1, \ldots, x_m) = 0$ is a semi-algebraic (resp. $k
[t]$-semi-algebraic)
condition, for $f$ a
polynomial over $k$ (resp. $k [t]$).

A subset of $\bar k ((t))^{m} \times \ZZ^{r}$
defined by a semi-algebraic (resp. $k [t]$-semi-algebraic)
condition is called {\emph {semi-algebraic}} (resp. {\emph
{$k [t]$-semi-algebraic}}).
One defines similarly semi-algebraic and $k [t]$-semi-algebraic
subsets of 
$K ((t))^{m} \times \ZZ^{r}$ for $K$ an algebraically closed field
containing
$\bar k$.

A function $\alpha : \bar k ((t))^m \times \ZZ^n \rightarrow \ZZ$
is called {\it simple} (resp. $k [t]$-{\it simple}) if its graph is semi-algebraic
(resp. $k [t]$-semi-algebraic).

We will use in an essential way the
following result on quantifier
elimination due to J. Pas \cite{Pas}.

\begin{theorem}\label{Pas}If $\theta$ is a semi-algebraic 
(resp. $k [t]$-semi-algebraic)
condition,
then $$(\exists \, x_1 \in \bar k ((t))) \,\, \theta (x_1, \ldots, x_m; \ell_1,
\ldots, \ell_r)$$
is semi-algebraic (resp. $k [t]$-semi-algebraic). Furthermore, for any
algebraically closed field $K$ 
containing
$\bar k$, 
$$(\exists \, x_1 \in K ((t))) \,\, \theta (x_1, \ldots, x_m; \ell_1,
\ldots, \ell_r)$$
is also semi-algebraic (resp. $k [t]$-semi-algebraic) 
and may be defined by the same conditions (\ie
independently of $K$).
\end{theorem}

\subsection{}Let $X$ be an
algebraic variety over $k$.
For $x \in \cL (X)$, we denote by $k_{x}$
the residue field of $x$ on $\cL (X)$,
and by $\tilde x$ the corresponding rational point
$\tilde x \in \cL (X) (k_{x}) = X (k_{x}[[t]])$.
When there is no danger of confusion we will often write $x$
instead of $\tilde x$.
A {\it semi-algebraic family of semi-algebraic subsets}
(resp. $k [t]$-{\it semi-algebraic family of $k [t]$-semi-algebraic subsets})
(for $n = 0$ a semi-algebraic subset (resp. $k [t]$-semi-algebraic subset))
$A_{\ell}$, $\ell \in \NN^n$, of $\cL (X)$ is a family of subsets
$A_{\ell}$ 
of $\cL (X)$ such that there exists a covering of $X$ by affine 
Zariski open
sets $U$ with
$$
A_{\ell} \cap \cL (U) =
\Bigl\{ x \in \cL (U) \bigm \vert
\theta (h_1 (\tilde x), \ldots, h_{m} (\tilde x); \ell)\Bigr\},
$$
where $h_1, \ldots, h_{m}$ are regular functions on 
$U$
and $\theta$ is a semi-algebraic condition
(resp. $k [t]$-semi-algebraic condition). 
Here the $h_{i}$'s and
$\theta$ may depend on $U$ and $h_{i} (\tilde x)$ belongs to $k_{x} [[t]]$.

Let $A$ be a semi-algebraic subset (resp. $k [t]$-semi-algebraic subset)
of $\cL (X)$.
A function $\alpha : A \times \ZZ^n
\rightarrow \ZZ \cup \{\infty \}$ is called {\it simple} (resp.
$k [t]$-{\it simple})
if the 
family of subsets $\{x \in \cL (X) \bigm | \alpha (x, \ell_1, \ldots,
\ell_n) = \ell_{n + 1}\}$,
$(\ell_1, \ldots,
\ell_{n + 1}) \in \NN^{n + 1}$,
is a semi-algebraic family of semi-algebraic subsets
(resp. a $k [t]$-{\it semi-algebraic family of $k [t]$-semi-algebraic subsets})
of $\cL (X)$.

We will use the following consequences of Theorem \ref{Pas}.

\begin{prop}\label{Pas2}
\begin{enumerate}
\item[(1)]If $X$ and $Y$ 
are algebraic varieties over $k$, $f : \cL (X) \rightarrow \cL (Y)$
is a $k [t]$-morphism
and $A$ is a $k [t]$-semi-algebraic subset of $\cL (X)$,
then $f (A)$ is a $k [t]$-semi-algebraic subset of $\cL (Y)$.
\item[(2)]
If $X$ is an algebraic variety over $k$
and $A$
is
a $k [t]$-semi-algebraic subset of $\cL (X)$, then $\pi_{n} (A)$ is a
constructible subset of $\cL_{n} (X)$.
\end{enumerate}
\end{prop}

\begin{proof}(1) is a direct consequence of Theorem \ref{Pas}.
The proof of (2) is similar to the proof of Proposition 2.3 in
\cite{Arcs}.
\end{proof}

\subsection{}By replacing $t$ by $t^{1 / d}$ in the definition, one
defines
similarly semi-algebraic (resp. $k [t]$-semi-algebraic)
subsets
of $\cL^{1 / d} (X)$.

\subsection{}\label{GK}We denote by $\cM$ the abelian group
generated by symbols $[S]$, for $S$ a variety over
$k$, with the relations $[S] = [S']$ if $S$ and $S'$ are
isomorphic and $[S] = [S'] + [S \setminus S']$
if $S'$ is closed in $S$. There is a natural ring structure
on $\cM$, the product being
induced by the cartesian product of varieties, and
to any constructible set $S$ in some variety  one naturally associates
a class $[S]$ in $\cM$.
We denote
by $\cM_{\rm loc}$ the localisation 
$\cM_{\rm loc} :=  \cM [\LL^{-1}]$ with
$\LL := [\AA^1_{k}]$. We denote by 
$F^m \cM_{\rm loc}$ the subgroup generated by
$[S] \LL^{- i}$ with ${\rm dim} \, S - i \leq -m$, and by
$\widehat{\cM}$ the completion
of $\cM_{\rm loc}$ with respect to
the filtration $F^{\cdot}$. We
will also denote by $F^{\cdot}$ the filtration induced on
$\widehat{\cM}$. 
We denote by $\overline \cM_{\rm loc}$ the image of 
$\cM_{\rm loc}$ in $\widehat \cM$.

\subsection{}In fact, for technical reasons appearing in the proof of
Lemma \ref{4.4}, we shall need
to consider the following
quotient 
$\cM_{/}$
of $\cM$, which is defined by adding the relation
$$[V / G] = [V],$$
for every vector space $V$ over $k$ endowed with a
linear action of a finite group $G$.
We shall still  denote by $\LL$ the class of the affine line, 
and, replacing $\cM$ by $\cM_{/}$, one defines similarly as above
rings
$\cM_{\rm loc /}$ and $\widehat \cM_{/}$.

\subsection{}Let $A$ be a $k [t]$-semi-algebraic subset of $\cL (X)$. We call
$A$ {\it weakly stable at level} $n \in \NN$ if $A$ is a union of fibers
of $\pi_{n} : \cL (X) \rightarrow \cL_{n} (X)$. We call $A$ {\it weakly
stable}
if it stable at some level $n$. Note that weakly stable
$k [t]$-semi-algebraic subsets form a boolean algebra.
Let $X$, $Y$ and $F$ be algebraic varieties over $k$,
and let
$A$, resp. $B$, be a constructible subset of $X$,
resp. $Y$. We say that
a map
$\pi : A \rightarrow B$ is a {\it piecewise morphism} if
there exists a finite partition of
the domain of $\pi$ into locally closed subvarieties of $X$
such that the restriction of $\pi$ to any of these
subvarieties is a morphism of
schemes.
We say that
a map
$\pi : A \rightarrow B$ is a
{\it piecewise trivial fibration with fiber}
$F$, if there exists a finite partition of $B$ in subsets $S$ which are
locally closed
in $Y$ such that $\pi^{- 1} (S)$ is locally closed in $X$ and
isomorphic, as
a variety over $k$, to $S \times F$, with $\pi$
corresponding under the isomorphism to the projection
$S \times F \rightarrow S$. We say that the map $\pi$ is
a
{\it piecewise trivial fibration over} some constructible subset $C$ of
$B$,
if the restriction of $\pi$ to $\pi^{- 1} (C)$ is a piecewise 
trivial fibration onto $C$.
One defines similarly the notion of a {\it pievewise vector bundle
of rank $e$}.

Let $X$ be an algebraic variety over $k$ of pure
dimension
$d$ (in particular we assume that
$X$ is non empty)
and let $A$ be a $k [t]$-semi-algebraic subset of $\cL (X)$. We call $A$
{\it stable at level} $n \in \NN$, if $A$ is weakly
stable at level $n$ and $\pi_{m + 1} (\cL (X)) \rightarrow \pi_{m} (\cL
(X))$
is a piecewise trivial fibration over $\pi_{m} (A)$ with fiber
$\AA^{d}_{k}$ for all $m \geq n$.
We call $A$ {\it stable} if it stable at some level $n$.

\begin{lem}\label{3.1}Let $X$ be an algebraic variety over $k$
of pure
dimension $d$, and let $A$ be a $k [t]$-semi-algebraic subset of $\cL (X)$. There
exists a reduced
closed subscheme 
$S$ of $X \otimes k[t]$, with
${\rm dim}_{k[t]} \, S < {\rm dim} \, X$,
and a $k [t]$-semi-algebraic family $A_i$, $i \in \NN$, of $k [t]$-semi-algebraic
subsets of $A$ such that
$\cL (S) \cap A$ and the $A_i$'s form a partition of
$A$, each $A_i$ is stable at some level $n_i$,
and
$$
\lim_{i \rightarrow \infty} \, ({\rm dim}\, \pi_{n_i} (A_i)
- (n_i + 1) \, d) = -\infty.
$$
Moreover, if $\alpha : \cL (X) \rightarrow \ZZ$ is a $k [t]$-simple function,
we can choose the partition such that
$\alpha$ is constant on
each $A_i$.
\end{lem}

\begin{proof}The proof is
literally the same as the one of Lemma 3.1 of \cite{Arcs}, noticing
that Lemma 4.4 of \cite{Arcs} also holds for a closed subscheme
$S$ of $X \otimes k [t]$
with ${\rm dim}_{k [t]} S < d$.
\end{proof}

Let $X$ be an algebraic variety over $k$ of pure
dimension $d$. Denote by $\BB^{t}$ the set of all
$k [t]$-semi-algebraic subsets of $\cL (X)$, and by
$\BB_{0}^{t}$ the set of all $A$ in $\BB^{t}$ which are stable. Clearly 
there is a unique additive measure
$$\tilde \mu : \BB_{0}^{t} \longrightarrow \cM_{\rm loc}$$
satisfying 
$$
\tilde \mu (A) = [\pi_{n} (A)] \, \LL^{- (n+1)d} \, ,
$$
when $A$ is stable at level $n$.

\begin{def-prop}
Let $X$ be an algebraic
variety
over $k$ of pure dimension $d$. Let $\BB^{t}$ be the set of all $k [t]$-semi-algebraic
subsets
of $\cL (X)$. There exists a unique map $\mu : \BB^{t} \rightarrow
\widehat \cM$ satisfying the following three properties.
\begin{enumerate}
\item[(1)] If $A \in
\BB^{t}$ is stable at level $n$, then
$\mu (A) = [\pi_{n} (A)] \LL^{- (n + 1) d}$.
\item[(2)] If $A \in
\BB^{t}$ is contained in $\cL (S)$ with $S$ a reduced
closed subscheme of $X \otimes k[t]$ with 
${\rm dim}_{k[t]} \, S < {\rm dim} \, X$, then $\mu (A) = 0$.
\item[(3)] Let $A_{i}$ be in $\BB^{t}$ for each $i$ in $\NN$.
Assume that the
$A_{i}$'s are mutually disjoint and that
$A := \bigcup_{i \in \NN} A_{i}$ is $k [t]$-semi-algebraic. Then
$\sum_{i \in \NN} \mu (A_{i})$ converges in $\widehat \cM$
to $\mu (A)$.
\end{enumerate}
Moreover we have:
\begin{enumerate}
\item[(4)]If $A$ and $B$ are in $\BB^{t}$, $A \subset B$ and if $\mu (B)
\in F^{m} \widehat \cM$, then 
$\mu (A)
\in F^{m} \widehat \cM$.
\end{enumerate}
This unique map $\mu$ is called the {\it motivic measure} on $\cL (X)$
and is denoted by $\mu_{\cL (X)}$ or $\mu$. 
For $A$ in $\BB^{t}$ and $\alpha : A \rightarrow \ZZ \cup
\{\infty\}$ a $k [t]$-simple
function,
one defines the motivic integral
$$
\int_{A} \LL^{- \alpha} d \mu := \sum_{n \in \ZZ} \mu (A \cap \alpha^{-1} (n))
\, \LL^{- n}
$$
in $\widehat \cM$, whenever the right hand side
converges in $\widehat \cM$, in which case we say that $\LL^{- \alpha}$
is integrable on $A$. If the function
$\alpha$ is bounded from below, then
$\LL^{- \alpha}$
is integrable on $A$, because of (4).
\end{def-prop}

\begin{proof}The proof of Definition-Proposition 3.2 of \cite{Arcs}
generalizes to the present case because Lemma 4.3
of \cite{Arcs} holds also for $X$
a scheme of finite type over $k
[t]$
(replacing ``dimension'' by ``relative dimension''),
and because Lemma 2.4 of \cite{Arcs} holds for
``semi-algebraic'' 
replaced by ``$k [t]$-semi-algebraic'' (cf. Lemma \ref{3} below), both
with identically
the same proofs. Note that we have to replace
Lemma 3.1 of \cite{Arcs} by Lemma \ref{3.1}.
\end{proof}

\subsection{}\label{omega}Let $h : \cL (Y) \rightarrow \cL (X )$ be a
$k [t]$-morphism
with $Y$ and $X$ of pure
dimension $d$.
Let $y$ be a closed point of $\cL (Y) \setminus
\cL (Y_{\rm sing})$ and denote by $\varphi$ the
corresponding
morphism
$\varphi : {\rm Spec} \, K [[t]] \rightarrow Y$, with $K$ a field
extension
of $k$. We define an element ${\rm ord}_t \cJ_{h} (y)$ in $\NN \cup
\{\infty\}$, the {\emph {order of the jacobian of $h$}} at $y$, as follows.
Consider the $K [[t]]$-module
$M = \varphi^{\ast} (\Omega^{d}_{Y})$ and set
$L := M \otimes_{K [[t]]} K ((t))$.  Here by $\Omega^{d}$ we mean $d$-th exterior power of the sheaf of
differentials.
The image $\tilde M$ of $M$ in the
$K ((t))$-vector space $L$
is a lattice of rank  1. One may also consider
the image $\tilde N$ of the module
$\varphi^{\ast} h^{\ast}(\Omega^{d}_{X \otimes k[t] | k [t]})$ in $L$.
If $\tilde N$ is non zero, $\tilde N = t^{n}\tilde M$ for some $n$ in
$\NN$ and one sets
${\rm ord}_t \cJ_{h} (y) = n$. When $\tilde N = 0$, one sets 
${\rm ord}_t \cJ_{h} (y) = \infty$.

Similarly assume  $Y$ is irreducible and let $\omega$ be
an element in $\Omega^{d}_{Y} \otimes_{k} k(Y)$. Denote by $\Lambda$ the
$K [[t]]$-submodule of $L$ generated by $\varphi^{\ast} (\omega)$.
If $\Lambda$ is non zero, $\Lambda = t^{n}\tilde M$ for some $n$ in
$\ZZ$ and one sets
${\rm ord}_t \omega (y) = n$. When $\Lambda = 0$, one sets 
${\rm ord}_t 
\omega (y) = \infty$.

\begin{lem}Let $X$ and $Y$ be 
$k$-varieties
of pure dimension $d$ and let
$h : \cL (Y) \rightarrow \cL (X )$
be a $k [t]$-morphism. Then the function $y \mapsto {\rm ord}_{t} \cJ_{h}
(y)$
is $k [t]$-simple
on $\cL (Y) \setminus
(\cL (Y_{\rm sing}))$. Similarly, if $Y$ is irreducible and
$\omega$ belongs to 
$\Omega^{d}_{Y} \otimes_{k} k(Y)$, the function
$y \mapsto {\rm ord}_{t} \omega (y)$
is $k [t]$-simple
on $\cL (Y) \setminus
\cL (Y_{\rm sing})$.

\end{lem}

\begin{proof}Direct verification.
\end{proof}

Under the preceeding assumptions, we extend the functions
${\rm ord}_t \cJ_{h} (y)$  and ${\rm ord}_t \omega (y)$  by
$\infty$
to
a
$k [t]$-simple function on $\cL (Y)$.

\begin{theorem}[Change of variables formula]\label{CV}
Let $X$ and $Y$ be algebraic varieties over $k$, of pure
dimension $d$. 
Let $h : \cL (Y) \rightarrow \cL (X )$ be a $k [t]$-morphism.
Let 
$A$ and $B$
be $k [t]$-semi-algebraic
subsets of $\cL (X)$ and $\cL (Y)$ respectively.
Assume that $h$ induces a bijection between
$B$ and $A$. Then,
for any $k [t]$-simple function
$\alpha : A \rightarrow
\ZZ \cup \{\infty\}$ such that $\LL^{- \alpha}$ is
integrable on $A$,
we have 
$$
\int_A \LL^{- \alpha} d \mu = 
\int_{B}
\LL^{- \alpha \circ h - {\rm ord}_t \cJ_{h} (y)} d \mu.
$$
\end{theorem}

\begin{proof}
By resolution of singularities we may assume that $Y$ is smooth. If 
$h$ is induced by a proper birational morphism from $Y$
to $X$, then Theorem \ref{CV} is a direct consequence of
Lemma 3.4 of \cite{Arcs}. In the general case it is a direct 
consequence of Lemma \ref{3.1}
and the following Lemma \ref{KL}.
\end{proof}

For $X$ a variety and $e$ in $\NN$,
we set
$$
\cL^{(e)} (X) := \cL (X) \setminus \pi_{e}^{-1} (\cL_{e} (X_{\rm sing})).
$$

We call a subset $A$ of $\cL (X)$ cylindrical at level $n$
if $A = \pi_{n}^{-1} (C)$, with $C$ a {\it constructible} subset of
$\cL_{n} (X)$. We say that $A$ is cylindrical if it is cylindrical
at some level $n$.

\begin{lem}\label{KL}Let $X$ and $Y$ be algebraic varieties over $k$, of pure
dimension $d$, with $Y$ smooth.
Let $h : \cL (Y) \rightarrow \cL (X )$ be a $k [t]$-morphism.
Let 
$B \subset \cL (Y)$
be cylindrical and put $A = h (B)$.
Assume that ${\rm ord}_t \cJ_{h} (\varphi)$
has constant value $e < \infty$ for all $\varphi \in B$, and that
$A \subset \cL^{(e')} (X)$ for some $e'$ in $\NN$.
Then $A$ is cylindrical. Morever,  if the restriction of $h$ to $B$
is injective,
then, for $n \in 
\NN$
large enough, we have the following:
\begin{enumerate}
\item[(a)]If $\varphi$ and $\varphi'$ belong to
$B$ and $\pi_{n} (h(\varphi)) = \pi_{n} (h(\varphi'))$,
then $\pi_{n - e}(\varphi) = \pi_{n - e}(\varphi')$.
\item[(b)]The morphism
$h_{n \ast} : \pi_{n} (B) \rightarrow \pi_{n} (A)$
induced by $h$ is a piecewise trivial fibration with fiber $\AA^e_{k}$.
\end {enumerate}
\end{lem}

\begin{proof}Let $n$ in $\NN$ be large enough.
We may assume that $B$ is
cylindrical at level 
$n - e$. That $A$ is cylindrical at level $n$
is an easy  consequence
of the following assertion:
\begin{enumerate}
\item[(a$''$)]For all $\varphi$ in $B$ and $x$ in $\cL (X)$, with
$\pi_{n} (h (\varphi)) = \pi_{n} (x)$, there exists $y$ in $\cL (Y)$
with $h (y) = x$ and $\pi_{n - e} (\varphi) =
\pi_{n - e} (y)$ (whence $y \in B$, since $B$
is cylindrical at level $n - e$).
\end{enumerate}

The proof of (a$''$) is the same as the proof of assertion (a$''$) in 
Lemma 3.4 of \cite{Arcs}.
(Note that with the notation of loc. cit. $B$ is contained in $\Delta_{e,
e'}$.)
Assertion (a) is a direct consequence of
(a$''$),
taking $x = h (\varphi')$ and using the injectivity of $h_{\vert 
B}$.
It remains to prove (b). Because of (a), we may assume that $X$ and $Y$ 
are affine.
Let $s : \cL_{n} (X) \rightarrow \cL (X)$ be a section
of the projection  
$\pi_{n} : \cL (X) \rightarrow \cL_{n} (X)$ such that the restriction 
of
$\pi_{n + e} \circ s$ to $\pi_{n} (A)$ is a piecewise morphism.
The existence of such a section has been shown in the proof of Lemma
3.4 of  \cite{Arcs}. Since $A$ is cylindrical at level $n$,
$s (\pi_{n} (A))$ is contained in $A$. Let $\theta$ be the mapping
$$
\theta : \pi_{n} (A) \longrightarrow B : \qquad x \longmapsto
h^{- 1} (s (x)).
$$
We will prove the following assertion:
\begin{enumerate}
\item[(c)]The map $\pi_{n} \circ \theta :
\pi_{n} (A) \rightarrow \pi_{n} (B)$ is a piecewise morphism.
\end{enumerate}

Using (c), the proof of (b) is the same as in the proof
of Lemma 3.4 in \cite{Arcs}, except that we have to replace the 
assertion that $\theta$ in loc. cit.
is a piecewise morphism by the slightly weaker assertion (c) above.

It only remains to prove (c). Let $x$ be in $\pi_{n} (A)$ and
$y$ in $\pi_{n} (B)$. Using assertion (a) we see that
$y = (\pi_{n} \circ \theta) (x)$ if and only if
there exists $\tilde y$ in $\pi_{n + e} (B)$ such that $y = \pi_{n} 
(\tilde y)$ and $h_{n + e \ast} (\tilde y) = \pi_{n + e} (s (x))$.
Thus, the graph of the map $\pi_{n} \circ \theta$ is constructible
and assertion (c) follows from the next lemma.
\end{proof}

\begin{lem}\label{cons}Let $X$ and $Y$ be algebraic varieties over $k$
and
let  $U$ and $V$ be constructible subsets of $X$ and $Y$ respectively.
If $f : U \rightarrow V$ is a map whose graph is a constructible 
subset of $X \times Y$, then $f$ is a piecewise morphism.
\end{lem}

\begin{proof}Well known.
\end{proof}

\begin{remark}All the material in this section (before \ref{cons})
generalizes
to ``$X$ and $Y$ algebraic varieties''
replaced by 
``$X$ and $Y$ separated reduced schemes of finite type over $k [t]$''.
In that case $X_{\rm sing}$ denotes the locus of points at which $X$ is
not smooth over $k [t]$, ``dimension'' has to be replaced by
``relative dimension over $k [t]$'', and in \ref{KL}
one replaces the hypothesis ``$Y$ smooth'' by
``$Y \otimes k(t)$ smooth''. Moreover one can also work with schemes
over $k [[t]]$ instead of over $k [t]$, replacing everywhere
$k [t]$ by $k [[t]]$. The proofs remain essentially the same, but
since this is not needed in the present paper, we do not give details here.
\end{remark}

\section{Study at the origin}\label{so}
\subsection{}\label{tet}Let $d \geq 1$ be an integer and
let $k$ be  field of characteristic 0 containing all $d$-th roots of unity.
Let $G$ be a finite subgroup of ${\rm GL}_{n} (k)$ of order $d$.
We fix a primitive $d$-th root of unity $\xi$
in $k$.
We denote by ${\rm Conj} (G)$ the set of conjugacy classes in $G$.
We let $G$ act on $\AA^{n}_{k}$ and we
consider
the morphism of schemes
$h : \tilde X = \AA^{n}_{k}
\longrightarrow
X = \AA^{n}_{k} / G.$
We denote by $ 0$ the origin in $\tilde X$ and $X$.
Let $\tilde
\Delta$ be the closed subvariety of $\tilde X$ consisting of the closed
points having a nontrivial stabilizer and let $\Delta$ be its image in
$X$
(the discriminant).
We denote by $\cL (X)^{g}$ (resp. $\cL^{1 / d} (\tilde X)^{g}$)
the complement of $\cL (\Delta)$ (resp.
$\cL^{1 / d} (\tilde \Delta)$)
in $\cL (X)$ (resp. $\cL^{1 / d} (\tilde X)$), and define similarly
$\cL (X)^{g}_{W}$ (resp. $\cL^{1 / d} ( \tilde X)^{g}_{W}$)
when $W$ is a subscheme of
$X$ (resp. $\tilde X$).

Let $\varphi$ be a geometric point of 
$\cL (X)^{g}_{0}$. So $\varphi$ is given by a morphism
$\varphi : {\rm Spec} \, K [[t]] \rightarrow X$ with $K$ an algebraically
closed overfield of $k$. The generic point of the image of $\varphi$ is
in $X \setminus \Delta$ and the special point is 0.
We can lift $\varphi$ to a morphism 
$\tilde \varphi$ making the following diagram commutative:
\begin{equation}\label{tildephi1}\xymatrix{
{\rm Spec} \, K [[t^{1 / d}]]
\ar[d] \ar[r]^<<<<{\tilde \varphi} & \tilde X \ar[d]^{h}\\
{\rm Spec} \, K [[t]]
\ar[r]^<<<<<{\varphi}&X.
}
\end{equation}

There is a unique
element $\gamma$ in $G$ such that
\begin{equation}\label{tildephi2}
\tilde \varphi (\xi t^{1 / d}) = \gamma \, \tilde \varphi (t^{1 / d}).
\end{equation}
If we change $\tilde \varphi$ in the diagram (\ref{tildephi1}),
$\gamma$ will be replaced by a conjugate.
If we denote by 
$\cL (X)^{g}_{0, \gamma}$ the set of $\varphi$'s in 
$\cL (X)^{g}_{0}$ such that there exists 
$\tilde \varphi$ satisfying (\ref{tildephi2}),
we have
$\cL (X)^{g}_{0, \gamma} = \cL (X)^{g}_{0, \gamma'}$
for $\gamma$ and $\gamma'$ in the same conjugacy class,
and we have a decomposition
$$\cL (X)^{g}_{0} =  \coprod \cL (X)^{g}_{0, \gamma} \, ,$$
for $\gamma$ running over a set of representatives of the conjugacy classes.

For each $\gamma$ in $G$, choose a basis $b_{\gamma}$ in which the
matrix
of $\gamma$ is diagonal, and denote by $\xi^{e_{\gamma, i}}$,
the diagonal coefficients,
with $1 \leq
e_{\gamma, i} \leq d$, $1 \leq i \leq n$.

\begin{lem}\label{lemme2.1.1}Let 
$\gamma$ be in $G$. A point $\tilde \varphi$
in $\cL^{1 / d} (\tilde X)^{g}$ projects to a point
in $\cL (X)^{g}_{0, \gamma}$ if and only if it is in the $G$-orbit
of a point in $\cL^{1 / d} (\tilde X)^{g}$ of
the form
\begin{equation}\label{2.1.3}
\tilde \varphi (t^{1 / d})=
(t^{e_{\gamma, 1} / d} \varphi_{1} (t),  \ldots, t^{e_{\gamma, n} / d} \varphi_{n} (t))
\end{equation}
in the basis $b_{\gamma}$.
\end{lem}

\begin{proof}It follows from (\ref{tildephi2}) that a point of
$\cL^{1 / d} (\tilde X)^{g}$ which projects to 
a point
in $\cL (X)^{g}_{0, \gamma}$ is in the $G$-orbit
of a point
of the form (\ref{2.1.3}).
To conclude observe that, in the basis $b_{\gamma}$,
$G$-invariant polynomials are sums of monomials of the form
$x_{1}^{m_{1}} \ldots x_{n}^{m_{n}}$, with $d$ dividing $\sum_{1 \leq i
\leq n} e_{\gamma, i} m_{i}$.
\end{proof}

\subsection{}We consider the morphism of $k [t]$-schemes
$$
\tilde \lambda : \AA^{n}_{k[t]}
\longrightarrow X \otimes k[t] 
\qquad
(x_{1}, \ldots, x_{n})
\longmapsto
h (t^{e_{\gamma, 1}/ d}x_{1},
\ldots,
t^{e_{\gamma, n}/ d}x_{n}),
$$ 
where $x_{1}, \ldots, x_{n}$ are the affine coordinates
corresponding to the basis $b_{\gamma}$. This is indeed
a $k [t]$-morphism, since, in the basis $b_{\gamma}$,
$G$-invariant polynomials are sums of monomials of the form
$x_{1}^{m_{1}} \cdots x_{n}^{m_{n}} $, with $d$ dividing
$\sum_{1 \leq i \leq n} e_{\gamma, i} m_{i}$.
The morphism $\tilde \lambda$ induces a $k [t]$-morphism
$\tilde \lambda_{\ast} : \cL (\AA^{n}_{k}) \rightarrow
\cL (X)_{0}$. Note that Lemma \ref{lemme2.1.1}
implies that
\begin{equation}\label{3.3.5}
\cL (X)^{g}_{0, \gamma}
=
\tilde \lambda_{\ast} (\cL (\AA^{n}_{k}))
\cap \cL (X)^{g}.
\end{equation}

\begin{prop}\label{3.4}For every $\gamma$ in $G$,
$\cL (X)^{g}_{0, \gamma}$ is a $k [t]$-semi-algebraic subset of $\cL (X)$.
\end{prop}

\begin{proof}This follows directly from
(\ref{3.3.5}) and Proposition
\ref{Pas2} (1).
\end{proof}

\subsection{}For $\gamma$ in $G$ we denote by $G_{\gamma}$
the centralizer of $\gamma$ in $G$. It follows from Theorem
\ref{Pas} that $\cL (\AA^{n}_{k}) / G_{\gamma}$ is a semi-algebraic
subset of
$\cL (\AA^{n}_{k} / G_{\gamma})$.

\begin{lem}\label{3.6}The morphism $\tilde \lambda$
is invariant under the action of $G_{\gamma}$ on $\AA^{n}_{k[t]}$.
Moreover the fibers of $\tilde \lambda_{\ast}$ above
$\cL (X)^{g}_{0, \gamma}$ are $G_{\gamma}$-orbits.
\end{lem}

\begin{proof}The first assertion is clear because the eigenspaces of
$\gamma$ are invariant subspaces under the action of $G_{\gamma}$.
Next we prove the second assertion. Let $x = (x_{1}, \ldots, x_{n})$ and
$x' = (x'_{1}, \ldots, x'_{n})$ be in a same fiber of
$\tilde \lambda_{\ast}$ above
$\cL (X)^{g}_{0, \gamma}$, and set 
$
\tilde \varphi =
(t^{e_{\gamma, 1} / d} x_{1},  \ldots, t^{e_{\gamma, n} / d}
x_{n})$ and 
$\tilde \varphi' =
(t^{e_{\gamma, 1} / d} x'_{1},  \ldots, t^{e_{\gamma, n} / d}
x'_{n})$.

Then (\ref{tildephi2}) holds for
$\tilde \varphi$, and also for $\tilde \varphi$
replaced by $\tilde \varphi'$.
There exists $\sigma$ in $G$ such that $\tilde \varphi' =
\sigma ( \tilde \varphi)$.
Hence (\ref{tildephi2}) also holds for $\tilde \varphi$
and $\gamma$ replaced by
$\sigma ( \tilde \varphi) = \tilde \varphi'$ and $\sigma \gamma \sigma^{-1}$
respectively.
Thus $\sigma = \sigma \gamma \sigma^{-1}$ and $\sigma \in G_{\gamma}$.
But the equality $\sigma ( \tilde \varphi) = \tilde \varphi'$ implies that
$\sigma (x) = x'$.
\end{proof}

By Lemma \ref{3.6}, $\tilde \lambda$ induces a morphism of $k
[t]$-schemes
$$
\lambda : (\AA^{n}_{k} / G_{\gamma}) \otimes
k [t] \longrightarrow X \otimes k [t].
$$
The morphism $\lambda$ induces a $k [t]$-morphism
$$
\lambda_{\ast} : \cL (\AA^{n}_{k} / G_{\gamma})
\longrightarrow \cL (X).
$$

\subsection{}\label{3.7}Considering 
$\cL (\AA^{n}_{k}) / G_{\gamma}$ as a (semi-algebraic)
subset of
$\cL (\AA^{n}_{k} / G_{\gamma})$ we have by (\ref{3.3.5})
and Lemma \ref{3.6} that $\lambda_{\ast}$ induces a bijection
between
$\cL (\AA^{n}_{k}) / G_{\gamma} \cap
\lambda_{\ast}^{-1} (\cL (X)^{g})$ and
$\cL (X)^{g}_{0, \gamma}$.

\section{Motivic Gorenstein measure of quotients}\label{mq}
\subsection{} Let $X$ be an irreducible
normal algebraic variety over $k$ of dimension
$d$
and assume
$X$ is Gorenstein with at most canonical singularities 
at each point. Hence there exists $\omega_{X}$ in
$\Omega^{d}_{X} \otimes_{k} k (X)$ generating
$\Omega^{d}_{X}$ at each smooth point of $X$, and, since $X$ is
canonical, ${\rm div} \, h^{\ast} (\omega_{X})$ is effective for any
resolution $h : Y \rightarrow X$. So by pulling back to $Y$ and using
the change of variables formula
Lemma 3.4 of \cite{Arcs}, we see that $\LL^{-\ord_{t} \omega_{X}}$
is integrable on $\cL (X)$ (see \ref{omega} for the definition of
$\ord_{t} \omega_{X}$). Furthermore
the function $\ord_{t} \omega_{X}$ does not depend
on the choice of $\omega_{X}$.
Hence one may define the
\textit{motivic Gorenstein measure}
$\mu^{Gor} (A)$ of a $k [t]$-semi-algebraic subset $A$ of $\cL (X)$ 
as
$$
\mu^{Gor} (A) :=
\int_{A} \LL^{- {\rm \ord}_{t} \omega_{X}} d \mu_{\cL (X)}
$$
in $\widehat \cM$.

\subsection{}\label{weight}Let $d \geq 1$ be an integer and
let $k$ be  field of characteristic 0 containing all $d$-th roots of unity.
Let $G$ be a finite subgroup of ${\rm SL}_{n} (k)$ of order $d$.
Set $X = \AA^{n}_{k} / G$ and let $h : \AA^{n}_{k} \rightarrow X$
be the projection. The variety
$X$ has only canonical Gorenstein singularities
and we can take
$\omega_{X}$ in $\Omega^{n}_{X / k} \otimes k (X)$ such that
$h^{\ast} (\omega_{X}) = dx_{1} \wedge \cdots \wedge dx_{n}$.

For $\gamma$ in $G$, recall the weight $w (\gamma)$ of
$\gamma$
by $w (\gamma) := \sum_{1 \leq i \leq n} e_{\gamma, i} / d$,
where the $e_{\gamma, i}$'s are as in \ref{tet}, \ie
$1 \leq e_{\gamma, i} \leq d$
and $\xi^{e_{\gamma, i}}$
are the eigenvalues of $\gamma$ for $ i = 1, \ldots, n$.
Note that $w (\gamma) \in \NN \setminus \{0\}$, since $G \subset
{\rm SL}_{n}(k)$.

\begin{lem}\label{4.3}For any $\gamma$ in $G$, we have
$$\mu^{Gor} (\cL (X)^{g}_{0, \gamma})
= \LL^{- w (\gamma)} \, 
\mu^{Gor}_{\cL (\AA^{n}_{k} / G_{\gamma})} (\cL (\AA^{n}_{k}) / G_{\gamma}).
$$
\end{lem}

\begin{proof}Let $\lambda$ be as in \S\kern .15em
\ref{so}. Direct verification yields
$\lambda^{\ast} (\omega_{X})
= t^{w (\gamma)} \omega_{\AA^{n}_{k} / G_{\gamma}}$.
The lemma follows now from \ref{3.7} and Theorem \ref{CV}
(with $h$ replaced by $\lambda_{\ast}$).
\end{proof}

A reason for considering the measure
$\mu_{\cL (X)}^{Gor}$ instead of
$\mu_{\cL (X)}$
is given by the next lemma. It is also at that place that it seems
necessary to work in the
ring
$\widehat \cM_{/}$ instead of just $\widehat \cM$.

\begin{lem}\label{4.4}The image of  $\mu_{\cL (X)}^{Gor}
(\cL (\AA^{n}_{k})/ G)$
in $\widehat \cM_{/}$ is equal to  $1$.
\end{lem}

\begin{proof}Let $M$ be a large integer.
For $e$ in $\NN$, we consider the subset
$\Delta_{e, M}$ of 
$\cL (\AA^{n}_{k})$ consisting of all
points $\varphi$ in 
$\cL (\AA^{n}_{k})$
such that ${\rm ord}_{t} \cJ_{h} (\varphi) = e$
and $h (\varphi) \in \cL^{(M)} (X)$.
Note that
$({\rm ord}_{t} \omega_{X}) \circ h = - {\rm ord}_{t} \cJ_{h}$,
because ${\rm ord}_{t} h^{\ast} (\omega_{X}) = {\rm ord}_{t}
(dx_{1} \wedge \cdots \wedge dx_{n}) = 0$.
Thus
$$
\mu_{\cL (X)}^{Gor}
(\cL (\AA^{n}_{k})/ G) = 
\sum_{e = 0}^{M} \,
\LL^{e}
\mu_{\cL (X)} (h (\Delta_{e, M})) + R_{M},
$$
with $\lim_{M \rightarrow \infty} R_{M} = 0$,
since $\LL^{-{\rm ord}_{t} \omega_{X}}$ is integrable on $\cL (X)$.
By the first assertion of Lemma \ref{KL} and by Lemma \ref{4.5} below,
for $m$ in $\NN$ large enough with respect to $M$,
we have for all $ e \leq M$ that $h (\Delta_{e, M})$
is stable at level $m$ and that 
$[\pi_{m}(h (\Delta_{e, M}))] = \LL^{- e}
[\pi_{m}(\Delta_{e, M}) / G]$.
Hence 
\begin{equation*}
\begin{split}
\mu_{\cL (X)}^{Gor}
(\cL (\AA^{n}_{k})/ G)  & = 
\sum_{e = 0}^{M} \,
[\pi_{m}(\Delta_{e, M}) / G] \, \LL^{- (m +1)n} + R_{M}
\\
& = [\pi_{m}(\cup_{e = 0, \ldots, M}
\Delta_{e, M}) / G] \, \LL^{- (m +1)n} + R_{M}
\\
& = [\pi_{m} (\cL(\AA^{n}_{k})) / G] \, \LL^{- (m +1)n} +
R'_{M},
\end{split}
\end{equation*}
with $\lim_{M \rightarrow \infty} R'_{M} = 0$
(because of Lemma 4.4 of \cite{Arcs}).
The lemma follows now, since
$\pi_{m} (\cL (\AA^{n}_{k}))/G$  is isomorphic to 
$\AA^{(m +1)n}_{k}/G$, the $G$-action on $\AA^{(m +1)n}_{k}$
being the diagonal one,
and the image of 
$\AA^{(m +1)n}_{k} / G$
in $\widehat \cM_{/}$ is equal to 
$\LL^{(m + 1)n}$ (it is here that we use the fact that we work in
$\widehat \cM_{/}$ instead of $\widehat \cM$).
\end{proof}

\begin{lem}\label{4.5}
Let $Y = \AA^{d}_{k}$ and 
$X = \AA^{d}_{k} / G$, with $G$ a finite subgroup of
${\rm GL}_{d} (k)$. Denote by 
$h : \cL (Y) \rightarrow \cL (X )$ the natural projection.
Let 
$B \subset \cL (Y)$
be cylindrical  and stable under the $G$-action.
Set $A = h (B)$.
Assume that ${\rm ord}_t \cJ_{h} (\varphi)$
has constant value $e < \infty$ for all $\varphi \in B$, and that
$A \subset \cL^{(e')} (X)$ for some $e'$ in $\NN$.
Then, for $n \in 
\NN$
large enough, we have the following:
\begin{enumerate}
\item[(a)]If $\varphi \in B$,
$\varphi' \in \cL (Y)$ and $\pi_{n} (h(\varphi)) = \pi_{n} (h(\varphi'))$,
then $\pi_{n - e}(\varphi)$ and $\pi_{n - e}(\varphi')$
have the same image in $\cL_{n - e}(Y) / G$.
\item[(b)]The morphism
$h_{n \ast} : \pi_{n} (B) / G \rightarrow
\pi_{n} (A)$ induced by $h$ may be endowed with the structure
of a piecewise
vector bundle of rank $e$.
\item[(c)]$[\pi_{n} (B) / G] = \LL^{e} \, [\pi_{n} (A)]$.
\end {enumerate}
\end{lem}

\begin{proof}Since assertion (a) is a direct consequence of assertion (a$''$)
in the proof of Lemma \ref{KL}, taking $x = h (\varphi')$, and
assertion (c) follows from (b),
it remains to prove (b).

By the first assertion in the statement of Lemma \ref{KL},
$A$ is cylindrical at level $n$, taking $n$ large enough.
In order to prove (b), we may assume
that $\pi_{n} (A)$ is a locally closed subvariey of $\cL_{n} (X)$.
The inverse image of $\pi_{n} (A)$ under the natural map $\cL_{n} (Y) / G
\rightarrow \cL_{n} (X)$ is locally closed
in $\cL_{n} (Y) / G$ and is equal to $\pi_{n} (B) / G$ by assertion (a)
and the fact that $B$ is cylindrical at level $n - e$, for $n$ large
enough. Hence $\pi_{n} (B) / G$ is a locally closed subvariety of
$\cL_{n} (Y) / G$, and $\pi_{n} (B)$ is a locally closed subvariety of
$\cL_{n} (Y)$.

Next we prove the following assertion:
\begin{enumerate}
\item[(d)]The stabilizer of $G$ acting on   $\pi_{n - e} (B)$
is trivial at every point of $\pi_{n - e} (B)$.
\end{enumerate}

Let $\sigma \in G \setminus \{1\}$
and set $\Delta_{\sigma} =\{y \in Y \, \vert \, \sigma (y) = y\}$.
Since $\ord_{t} \cJ_{h} \not= \infty$
on $B$, we have $B \cap \cL (\Delta_{\sigma}) = \emptyset$. Hence $B$ is
contained in
$\cup_{m \in \NN} (\cL (Y) \setminus \pi_{m}^{-1}
(\cL_{m} (\Delta_{\sigma})))$.
Thus, since $B$ is cylindrical, Lemma \ref{3} implies that $B$ is
contained in 
$\cL (Y) \setminus \pi_{m}^{-1}
(\cL_{m} (\Delta_{\sigma}))$
when $m $ is large enough. This concludes the proof of assertion (d).

Our next step is to construct a section of
the morphism
$h_{n \ast} : \pi_{n} (B) / G \rightarrow
\pi_{n} (A)$.
Let $s : \cL_{n} (X) \rightarrow \cL (X)$ be a section of the projection
$\pi_{n} : \cL (X) \rightarrow \cL_{n} (X)$ such that the restriction of
$\pi_{n + e} \circ s$ to $\pi_{n} (A)$ is a piecewise morphism. The
existence of such a section $s$ has been shown in the proof of Lemma
\ref{KL}.
Note that $s
(\pi_{n} (A))$ is contained in $A$, since $\pi_{n} (A)$ is cylindrical
at level $n$. Denote by 
$\theta$ the map
$$
\theta : \, \pi_{n} (A) \longrightarrow B / G
\, : 
\quad x \longmapsto h^{-1} (s (x)) \mod G,
$$
and set
$$
\overline \theta = \tilde \pi_{n} \circ \theta : \,
\pi_{n} (A) \longrightarrow \pi_{n} (B) / G
\, : 
\quad x \longmapsto \theta (x) \mod t^{n + 1},
$$
where $\tilde \pi_{n} : B/ G \rightarrow \pi_{n} (B) / G$ is the
projection.
Clearly $\overline \theta$ is a section of $h_{n \ast}$.
One proves that $\overline \theta$ is a piecewise morphism by exactly
the same argument as for assertion (c) in the proof of Lemma \ref{KL},
replacing $B$, $\pi_{n} (B)$, and $\pi_{n + e} (B)$ by their
quotient under the action of $G$.

By (d), the natural morphism $p : \pi_{n} (B) \rightarrow \pi_{n} (B) /
G $
is {\'e}tale.
We consider the fiber product
$$\widetilde{\pi_n (A)} := 
\pi_n (A) \times_{\pi_n (B) /G}
\pi_n (B).$$
The strategy of proof is to construct a $G$-equivariant morphism
$\gamma : \pi_n (B) \rightarrow \widetilde{\pi_n (A)}$,
such that  the following diagram is commutative,
\begin{equation}\xymatrix{
\pi_n (B) / G
\ar@<0.5ex>[r]^>>>>>{h_{n \ast}} & \ar@<0.5ex>[l]^>>>>>{\overline \theta}
\pi_n (A)\\
\pi_n (B) \ar[u]^p
\ar@{.>}@<0.5ex>[r]^<<<<<<<{\gamma}& \ar@<0.5ex>[l] \widetilde{\pi_n (A)},\ar[u]
}
\end{equation}
then
to show it may  be endowed with  the structure
of a piecewise
vector bundle of rank $e$, and finally
to conclude by {\'e}tale descent.

We first construct the mapping $\gamma$.
Let  $\varphi$ be a point in $\pi_{n} (B)$.
It follows from (a) that there exists
a lifting
$\widetilde \varphi$  in $\pi_{n} (B)$
of 
$\overline \theta (h_{n \ast} (p (\varphi)))$
such that
$\varphi \equiv \widetilde \varphi \mod t^{n + 1 -e}$.
Furthermore, by (d), the
lifting
$\widetilde \varphi$ is uniquely determined by
$\varphi$. We set
$$
\gamma (\varphi) := 
(h_{n \ast} (p (\varphi)), \widetilde \varphi).
$$
Clearly, the graph of $\gamma$ is
constructible,
hence, by Lemma \ref{cons},
$\gamma$ is a piecewise morphism. We shall show later that,
as soon as $\overline \theta$ is a morphism and $\pi_n (B)$ is smooth,
$\gamma$ is
an actual  morphism.
Now take a point $(a, \widetilde \varphi)$ in $\widetilde{\pi_n (A)}$. 
We have $a = h_{n \ast} (p (\widetilde \varphi))$ and 
$\widetilde \varphi$ is a lifting of 
$\overline \theta (a)$.
Hence the conditions for a point $\varphi$ to be in the fiber
$\gamma^{-1}(a, \widetilde \varphi)$ are that
$\varphi \equiv \widetilde \varphi \mod t^{n + 1-e}$
and
$h (\varphi) \equiv h (\widetilde \varphi) \mod t^{n + 1}$.
Rewriting the first condition as
$\varphi = \widetilde \varphi + t^{n + 1-e} u$, with a unique
$u$ in $\cL_{e - 1} (\AA^d_k)$, the fiber
$\gamma^{-1}(a, \widetilde \varphi)$
can be determined by rewriting
the condition
$$
h (\widetilde \varphi + t^{n + 1 - e}u) \equiv h (\widetilde \varphi) \mod t^{n + 1}
$$
using the Taylor expansion of $h$ at $\widetilde \varphi$.
In this way, using again that
$n$ is large enough and that
$B$ is cylindrical at level $n - e$,
we find that
$$
\gamma^{-1}(a, \widetilde \varphi)
=
\Bigl\{\widetilde \varphi + t^{n + 1 - e}
(u_{0} + u_{1}t + \cdots + u_{e - 1}t^{e - 1})
\Bigm | 
L_{\widetilde \varphi} (u_{0}, \ldots, u_{e - 1}) = 0
\Bigr\},
$$
where $L_{\widetilde \varphi} (u_{0}, \ldots, u_{e - 1}) = 0$ is 
a system of linear homogeneous
equations
whose coefficients are regular functions of $\widetilde \varphi \in \cL_{n} (Y)$.

We refer to \cite{Arcs} 3.4 (3) for more details.
Moreover the solution space of this linear system has dimension $e$, since
the jacobian
matrix of $h$ at any point in 
$\pi^{-1}_{n} (\widetilde \varphi) $
is equivalent over $\bar k [[t]]$ to a diagonal matrix with diagonal
elements
$t^{e_{1}}$, $t^{e_{2}}$, \dots, satisfying $e = e_{1} + e_{2} +
\cdots$,
cf. \cite{Arcs} 3.4 (4).

In order to prove (b),
we may assume 
that $\pi_{n} (A)$ is a locally closed smooth 
subvariety of $\cL_{n} (X)$ and that $\overline \theta$ is a morphism,
provided 
that from now
on we only assume $B$ is cylindrical at level $n$ and that
we do not anymore increase $n$, which could destroy
the property of $\overline \theta$ to be a morphism.
When $k = \CC$, we see from our previous discussion about
$\gamma^{-1}(a, \widetilde \varphi)$, that
$\pi_{n} (B)$ is locally bianalytically isomorphic to
$\pi_n (A) \times \CC^{e}$. Hence $\pi_{n} (B)$
is smooth  for any $k$.
Now let us prove that $\gamma$ is a morphism.
When $k = \CC$, it is easy to see that $\gamma$
is continuous, hence is a morphism, since its domain is smooth and
it is a piecewise morphism. 
Thus
by the Lefschetz
principle, it follows that $\gamma$ is a morphism, for any $k$.
The fact that it may be endowed with the structure of a vector bundle
of rank $e$ follows from the above description of the fibers.
Now by {\'e}tale descent (Hilbert's Theorem 90, see, {\it e.g.},
\cite{Milne} p.124), we deduce that $h_{n \ast}$ may be endowed
with the structure of a 
vector bundle
of rank $e$.
\end{proof}

\bigskip

We can now prove the main result.

\begin{theorem}\label{MT}Let $d \geq 1$ be an integer and
let $k$ be  field of characteristic 0 containing all $d$-th roots of unity.
Let $G$ be a finite subgroup of ${\rm SL}_{n} (k)$ of order $d$,
so
$G$ acts on $\AA^{n}_{k}$. Consider the quotient $X := \AA^{n}_{k} /
G$.
\begin{enumerate}
\item[(1)]For any $\gamma$ in $G$, we have
$$\mu^{Gor} (\cL (X)^{g}_{0, \gamma})
= \LL^{- w (\gamma)}
$$
in
$\widehat \cM_{/}$.
\item[(2)]The relation 
\begin{equation*}\label{ME}
\mu^{Gor} (\cL (X)_{0})
= \sum_{[\gamma] \in {\rm Conj} (G)} \, \LL^{- w (\gamma)},
\end{equation*}
holds in the ring
$\widehat
\cM_{/}$,
where 
${\rm Conj} (G)$ denotes the set of conjugacy classes in $G$.
\end{enumerate}
\end{theorem}

\begin{proof}The first statement
is a direct consequence of Lemma \ref{4.3} and Lemma
\ref{4.4}
with $G$, $X$ replaced by $G_{\gamma}$, $\AA^{n}_{k} / G_{\gamma}$.
The second follows then, using
the
decomposition
$$\cL (X)^{g}_{0} =  \coprod \cL (X)^{g}_{0, \gamma}$$
and the fact that $\mu^{Gor} (\cL (X)_{0} \setminus \cL (X)^{g}_{0}) = 0$.
\end{proof}

\subsection{}\label{orb}Keeping the above notations,
we now assume 
that $G$ is a finite subgroup of ${\rm GL}_{n} (k)$, instead of 
${\rm SL}_{n} (k)$. Notice that now the weight $w (\gamma) \in
\QQ$ of an element $\gamma$ in $G$ might not be integral and that
$\omega_{X}$ might not exist. To remedy this  we consider the function
$\alpha_{X} : \cL (X) \rightarrow \QQ \cup \{\infty\}$
which is
defined by $\alpha_{X} (\varphi) = - \ord_{t} \cJ_{h} (\tilde \varphi)$
for any $\tilde \varphi$ in $\cL^{1/d} (\AA^{n}_{k})$ with $h (\tilde
\varphi) = \varphi$. Clearly $\alpha_{X} = \ord_{t} \omega_{X}$
when $G \subset {\rm SL}_{n} (k)$.
We define the {\emph {motivic orbifold measure}}
$\mu^{orb} (A)$ of a $k [t]$-semi-algebraic subset $A$ of $\cL (X)$ as
$$
\mu^{orb} (A) := \int_{A} \LL^{- \alpha_{X}} \, d\mu_{\cL (X)} \,
\in \widehat \cM \, [\LL^{1/d}].
$$
Theorem \ref{MT} remains true for
$G \subset {\rm GL}_{n} (k)$ if we replace $\mu^{Gor}$ by $\mu^{orb}$
and $\widehat
\cM_{/}$ by $\widehat
\cM_{/} \, [\LL^{1/d}]$. Indeed the proofs remain
basically the same, replacing $\ord_{t} \omega_{X}$ by $\alpha_{X}$.
At
the
same time one verifies that the integrals 
$\mu^{orb} (A)$ converge.
At the level of Hodge realization a
similar result is contained in \S\kern .15em 7 of \cite{Ba}.
Indeed, with the notation of loc. cit.,
$H (\mu^{orb} (\cL (X))) = E_{st}
(X, \Delta_{X}; u, v)$.

\subsection{}\label{add}More generally we may consider
a smooth irreducible
algebraic variety
$\tilde X$ endowed with an effective action
of a finite group $G$ of order $d$. We assume the field $k$
contains all $d$-th roots of unity.
We shall also assume that every $G$-orbit is
contained
in an affine open subset of $\tilde X$ and we
denote by $X$ the quotient variety $\tilde X / G$.
Using the previous methods, it is possible to express
$\mu^{orb} (\cL (X))$ in terms of weights associated to the group
action along the orbifold strata, similarly as in
\S\kern .15em 7 of \cite{Ba}, cf. \cite{Loo}.

\section{Chow motives and realizations}\label{mot}

\subsection{}We denote by
$\cV_{k}$ the category of smooth and projective $k$-schemes.
For  an object $X$ in $\cV_{k}$ and an integer $d$,
we denote by  $A^{d} (X)$ the Chow group of codimension
$d$ cycles with rational coefficients
modulo rational equivalence.
Objects of the category ${\rm CHM}_{k}$ of (rational)
$k$-motives
are triples $(X, p, n)$ where $X$ is in $\cV_{k}$,
$p$ is an idempotent (\ie $p^{2} = p$) in the ring of
correspondences ${\rm Corr}^{0} (X, X)$
($= A^{d} (X \times X)$ when $X$ is of pure dimension $d$), and
$n$ is an integer. If $(X, p, n)$
and $(Y, q, m)$ are motives, then
$$
{\rm Hom}_{{\rm CHM}_{k}} ((X, p, n), (Y, q, m))
=
q \, {\rm Corr}^{m - n} (X, Y) \, p.
$$
Here ${\rm Corr}^{r} (X, Y)$ is the group of correspondences of degree
$r$ from $X$ to $Y$ (which is $A^{d +r} (X \times Y)$ when $X$ is of
pure dimension $d$).
Composition of morphisms is given by composition of correspondences.
The category ${\rm CHM}_{k}$ is  additive, $\QQ$-linear, and pseudo-abelian,
and 
there is a natural tensor product on ${\rm CHM}_{k}$.
We denote by $h$ the functor $h : \cV_{k}^{\circ} \rightarrow
{\rm CHM}_{k}$ which sends an object $X$ to $h (X) = (X, {\rm id}, 0)$
and a morphism $f : Y \rightarrow X$ to its graph in
${\rm Corr}^{0} (X, Y)$. We denote by ${\LL}$ the Lefschetz motive
$\LL = ({\rm Spec} \, k, {\rm id}, -1)$. There is a canonical isomorphism
$h (\PP^{1}_{k}) \simeq 1 \oplus \LL$. 

Let $K_{0} ({\rm CHM}_{k})$ be the Grothendieck group of the pseudo-abelian
category
${\rm CHM}_{k}$. It is also the abelian group associated to the monoid of
isomorphism classes of $k$-motives
with respect to the addition $\oplus$. The tensor product on ${\rm CHM}_{k}$ induces 
a natural ring structure on $K_{0} ({\rm CHM}_{k})$.
For $m$ in $\ZZ$, let $F^{m} K_{0} ({\rm CHM}_{k})$
denote 
the subgroup of $K_{0} ({\rm CHM}_{k})$ generated by
$h (S, f, i)$, with $i - \dim S \geq m$. This gives a filtration
of the ring $K_{0} ({\rm CHM}_{k})$
and we denote by 
$\widehat K_{0} ({\rm CHM}_{k})$
the completion of 
$K_{0} ({\rm CHM}_{k})$ with respect to this filtration.

Gillet and Soul{\'e} \cite{G-S} and
Guill{\'e}n and Navarro Aznar \cite{G-N} proved the following result.
\begin{theorem}Let $k$ be a field
of characteristic 0.
There exists a unique map
$\chi_{c}$ which to any variety $X$ over $k$ associates
$\chi_{c} (X)$ in 
$K_{0} ({\rm CHM}_{k})$
such that
\begin{enumerate}
\item[(1)] If $X$ is smooth and projective,
$\chi_{c} (X) = [h (X)]$.
\item[(2)] If $Y$ is a closed reduced
subscheme in a variety
$X$
$$\chi_{c} (X \setminus Y) = \chi_{c} (X) - \chi_{c} (Y).$$
\item[(3)] If $X$ is 
a variety, $U$ and $V$ are open reduced
subschemes of $X$,
$$
\chi_c (U \cup V)
=
\chi_c (U) +
\chi_c (V) -
\chi_c (U \cap V).
$$
\item[(4)]If $X$ and $Y$ are varieties
$$
\chi_{c} (X \times Y) = \chi_{c} (X)
\,
\chi_{c} (Y).$$
\end{enumerate}
Furthermore, $\chi_{c}$ is determined by conditions
(1)-(2).
\end{theorem}
Hence $\chi_{c}$ induces a morphism of rings
$\chi_{c} : \cM \rightarrow K_{0} ({\rm CHM}_{k})$ with
$\chi_{c} (\LL) = \LL$
and
extends to a morphism
$\widehat
\chi_{c} : \widehat \cM \rightarrow \widehat
K_{0} ({\rm CHM}_{k})$.

\subsection{}\label{hr}
Recall that the Hodge polynomial of an algebraic variety $S$ defined
over a subfield of $\CC$
is the polynomial
$$
H (S; u, v) := \sum_{p, q} e^{p, q} (S) \, u^p v^q
$$
with
$$
e^{p, q} (S) := \sum_{i \geq 0} (- 1)^i h^{p, q} (H^i_c (S, \CC)),
$$
where $h^{p, q} (H^i_c (S, \CC))$ denotes the rank of the
$(p, q)$-Hodge component of the $i$-th cohomology group
with compact supports. One defines similarly the Hodge polynomial of
Chow motives.
It follows from a weight argument,
cf. \cite{Arcs} and \cite{Motivic}, that
the Hodge polynomial $H$
factorizes (hence also the 
Euler Characteristic $\rm Eu$)
through the image of  
$K_0 ({\rm CHM}_{k})$
in
$\widehat K_0 ({\rm CHM}_{k})$.

\subsection{}The following proposition shows that the morphisms
$\chi_c$ and $\widehat \chi_c$
factorize
through $\cM_{/}$ and
$\widehat \cM_{/}$ respectively.

\begin{prop}\label{quot}Let $V$ be a finite
dimensional vector space over $k$ and let $G$
be a finite subgroup of
${\rm GL} (V)$. Then the  following equality holds:
$$
\chi_{c} (V / G) = \chi_{c} (V).
$$
\end{prop}

\begin{proof}We will use the functor $h_{c}$
of
\cite{G-N} which to a variety $X$ over $k$ associates
an object $h_{c} (X) $ of the homotopy category ${\rm Ho} (C^{b}
({\rm CHM}_{k}))$ of bounded complexes
of objects in ${\rm CHM}_{k}$, such that $\chi_{c} (X)$ is the Euler
characteristic of
$h_{c} (X) $. Consider the functor $\tau : {\rm CHM}_{k} \rightarrow
{\rm Ho} (C^{b}
({\rm CHM}_{k})) $ which to an object $M$ associates 
the complex
in ${\rm Ho} (C^{b}
({\rm CHM}_{k}))$ which is zero
in non zero degree and is equal to $M$ in degree 0.
It follows from the identity 
$h (\PP^{1}_{k}) \simeq 1 \oplus \LL$ in ${\rm CHM}_{k}$ and the
definitions
that $h_{c} (V) $ is isomorphic to $\tau (\LL^{{\rm dim} V})$
in
${\rm Ho} (C^{b}
({\rm CHM}_{k}))$. By Corollary 5.3 of \cite{B-N}, 
$h_{c} (V / G) $ is a direct factor of $h_{c} (V) $
in ${\rm Ho} (C^{b}
({\rm CHM}_{k}))$. The functor $\tau$ being fully faithful 
and $\LL^{r}$ being indecomposable, it follows that
$h_{c} (V / G) $ is zero or equal to $\tau (\LL^{{\rm dim} V})$. 
Using a realization, for instance the Betti realization, 
one obtains that $h_{c} (V / G) = \tau (\LL^{{\rm dim} V}) $, and the
result follows.
\end{proof}

\section{Relation with resolution of singularities and
the McKay correspondence}\label{mk}

Let $X$ be an algebraic variety over $k$ of pure
dimension
$d$, and let 
$h : Y \rightarrow X$ be a resolution of singularities of
$X$. By this we
mean
$Y$ is a smooth algebraic variety over $k$, $h$ is birational,
proper and defined over $k$, and the exceptional
locus $E$ of $h$ has normal crossings, meaning that the $k$-irreducible
components of $E$ are smooth and intersect transversally.
Let us denote the $k$-irreducible components of $E$ by $E_{i}$, $i \in
J$.
For $I \subset J$, set $E_{I} = \bigcap_{i \in I} E_{i}$ and
$E_{I}^{\circ} = E_{I} \setminus \bigcup_{j \not\in I} E_{j}$.
Assume now $X$ is 
Gorenstein with at most canonical singularities 
at each point and consider $\omega_{X}$ in
$\Omega^{d}_{X} \otimes_{k} k (X)$ generating
$\Omega^{d}_{X}$ at each smooth point of $X$.
For $i$ in $I$, we
denote by $\nu_i  - 1$ the length of
$\Omega^d_{Y} / h^{\ast} \omega_{X} \cO_{Y}$ at the generic point
of $E_{i}$.

Let $W$ be a closed subvariety of $X$. By
Lemma 3.3 of \cite{Arcs} (cf. Proposition 6.3.2 of
\cite{Arcs}), the following
formula holds in $\widehat \cM$:
\begin{equation}\label{ast}
\mu^{Gor} (\pi_{0}^{-1} (W))
=
\LL^{- d} \sum_{I \subset J} \, [E_{I}^{\circ} \cap h^{-1} (W)] \,
\prod_{i \in I} \frac{\LL - 1}{\LL^{\nu_{i}} -1} \, .\tag{$*$}
\end{equation}

Now we can specialize to the case
where $X = \AA^{n}_{k} /G$ with $G$ a finite subgroup
of ${\rm SL}_{n} (k)$ and $W = \{0\}$. Theorem \ref{MT}
may now be rephrazed as follows.

\begin{theorem}\label{reph}Let $d \geq 1$ be an integer and
let $k$ be  field of characteristic 0 containing all $d$-th roots of unity.
Let $G$ be a finite subgroup of ${\rm SL}_{n} (k)$ of order $d$.
Let $h : Y \rightarrow X$ be a resolution of $X = \AA^{n}_{k} / G$.
Then the following relation holds in
$\widehat \cM_{/}$:
$$\LL^{- n} \sum_{I \subset J} \, 
[E_{I}^{\circ} \cap h^{-1} (0)] \,
\prod_{i \in I} \frac{\LL - 1}{\LL^{\nu_{i}} -1}
=
\sum_{[\gamma] \in {\rm Conj} (G)} \, \LL^{- w (\gamma)}. \hfill \qed$$
\end{theorem}

In particular, if the resolution $h$ is crepant, \ie all the $\nu_{i}$'s
are equal to 1, we get as a corollary the following form of the
McKay correspondence (cf. \cite{Reid}).

\begin{cor}Let $h : Y \rightarrow X$ be a crepant
resolution of $X = \AA^{n}_{k} / G$.
Then the following relation holds in
$\widehat \cM_{/}$:
$$[h^{-1} (0)] \,
=
\sum_{[\gamma] \in {\rm Conj} (G)} \, \LL^{n - w (\gamma)}. \hfill \qed$$
\end{cor}

By passing to the Hodge realization, cf. \ref{hr},
one obtains in particular the 
following  form of the
McKay correspondence, which was conjectured by Reid in \cite{Reid}
and proved by Batyrev in \cite{Ba}, see also \cite{Ba4}, \cite{Re}.

\begin{cor}\label{LC}Let $h : Y \rightarrow X$ be a crepant
resolution of $X = \AA^{n}_{k} / G$.
Then 
$$
H (h^{-1} (0)) \,
=
\sum_{[\gamma] \in {\rm Conj} (G)} \, (uv)^{n - w (\gamma)} \quad
\text{and}
\quad
{\rm Eu} (h^{-1} (0))) = {\rm card} \, {\rm Conj} (G). \hfill
\qed
$$
\end{cor}

\begin{remark}In the framework of \ref{orb},
when $A$ is of the form $\pi_0^{-1} (W)$,
one may express $\mu^{orb} (A)$ in terms of a resolution
of $X$ in a way completely similar to 
(\ref{ast}), replacing the integers $\nu_i$ 
by rational numbers
$\nu_i^{\ast}$ similarly defined with the help of  $\alpha_X$, cf. \cite{Loo}.
\end{remark}

\appendix

\setcounter{equation}{0}

\section{Measurable subsets of $\cL (X)$}

Let $X$ be an algebraic variety of pure dimension $d$
over a field $k$ of characteristic zero. We develop 
here the theory of measurable subsets of $\cL (X)$. When $X$ is smooth,
a measure theory for $\cL (X)$ in the case of the Hodge realization
has been considered  by Batyrev in
\cite{Ba2}.

\subsection{}We call a cylindrical subset $A$ of $\cL (X)$
stable at level $n \in \NN$ if $A$ is cylindrical
at level $n$ and $\pi_{m + 1} (\cL (X)) \rightarrow \pi_{m} (\cL (X))$
is a piecewise
trivial fibration over $\pi_{m} (A)$ with fiber $\AA^{d}_{k}$ for all $m
\geq n$.
We call $A$ stable if it is stable at some level $n$.

Denote by $\CC_{0}$ the family of stable cylindrical subsets of $\cL
(X)$
and by $\CC$ the boolean algebra  of cylindrical subsets of $\cL (X)$.
Since there might exist cylindrical subsets of $\cL (X)$ which are not
semi-algebraic, we cannot apply the motivic measure $\mu$ of   
\S\kern .15em \ref{pre} to elements of $\CC$ or $\CC_{0}$.
Some precautions are necessary.

\subsection{}Clearly there exists a unique additive measure
$$
\tilde \mu : \CC_{0} \longrightarrow \cM_{\rm loc}
$$
satisfying
$$
\tilde \mu (A) = [\pi_{} (A)] \, \LL^{- (n + 1) d}
$$
when $A \in \CC_{0}$ is stable at level $n$.
For $A$ in $\CC$, we define
$$
\mu (A) = \lim_{e \rightarrow \infty} \tilde \mu (A \cap \cL^{(e)} (X))
\in \widehat \cM.
$$
Indeed,
$A \cap \cL^{(e)} (X)$ is stable by Lemma 4.1 of \cite{Arcs},
and the limit exists in $\widehat \cM$ by Lemma 4.4 of loc. cit.
Moreover if $A \in \CC_{0}$ then $\mu (A)$ is the image in $\widehat
\cM$ of $\tilde \mu (A)$. Clearly $\mu$ is additive on $\CC$,
and even $\sigma$-additive because of the following lemma,
which first appeared in \cite{Arcs} Lemma 2.4 for weakly
stable semi-algebraic subsets, with a proof which actually holds also
for cylindrical
subsets. A different proof is given in Theorem 6.6 of \cite{Ba2}.

\begin{lem}\label{3}Let $A_{i}$, $i \in \NN$, be a family of cylindrical
subsets of $\cL (X)$. Suppose that $A := \cup_{i \in \NN} A_{i}$
is cylindrical. Then $A$ equals the union of a finite number of the $A_{i}$'s.
\end{lem}

\subsection{}We consider on $\widehat \cM$ the norm $||\cdot||$
defined by
$$
||\cdot|| : \quad \widehat \cM \longrightarrow \RR_{\geq 0} \quad :
\qquad 
a \longmapsto ||a|| := 2^{-n},
$$
where $n$ is the largest $n$ such that $a \in F^{n} \widehat \cM$.

For all $a$, $b$ in $\widehat \cM$, we have
$||ab|| \leq ||a|| \, ||b||$ and 
$||a + b|| \leq \max \, (||a||, ||b||)$.

Note also that, for all $A$, $B$ in $\CC$,
we have
$$
||\mu (A \cup B) || \leq \max \, (||\mu (A)||, ||\mu (B)||)
$$
and $||\mu (A)|| \leq ||\mu (B)||$ when $A \subset B$.

For $A$ and $B$ subsets of the same set, we use the notation
$A \triangle B$ for $A \cup B \setminus A \cap B$.

\begin{definition}We say that a subset $A$ of $\cL (X)$ is
{\emph{measurable}}
if, for every positive real number $\varepsilon$, there exists a sequence of
cylindrical subsets
$A_{i} (\varepsilon)$, $i \in \NN$,
such that
$$
\Bigl(A \triangle
A_{0} (\varepsilon)
\Bigr)
\subset \bigcup_{i \geq 1} A_{i} (\varepsilon),
$$
and 
$||\mu (A_{i} (\varepsilon))|| \leq \varepsilon$
for all $i \geq 1$. We say that $A$ is {\emph{strongly measurable}}
if moreover we can take $A_{0} (\varepsilon) \subset A$.
\end{definition}

\begin{theorem}If $A$ is a measurable subset of $\cL (X)$,
then 
$$\mu (A) :=
\lim_{\varepsilon \rightarrow 0}  \mu (A_{0} (\varepsilon))$$
exists in $\widehat \cM$ and is independent of the choice of the
sequences
$A_{i} (\varepsilon)$, $i \in \NN$.
\end{theorem}

\begin{proof}This is proved in exactly the same way as Theorem 6.18 of
\cite{Ba2} using Lemma \ref{3}.
\end{proof}

For $A$ a measurable subset of $\cL (X)$, we shall call
$\mu (A)$ the motivic measure of $A$.

One should remark that obviously
any cylindrical subset of $\cL (X)$ is strongly measurable and that
the measurable subsets of $\cL (X)$
form a boolean algebra.
Note also that if $A_{i}$, $i \in \NN$, is a sequence
of measurable subsets
of $\cL (X)$ with
$\lim_{i \rightarrow \infty} ||\mu (A_{i})|| = 0$, then
$\cup_{i \in \NN} A_{i}$ is measurable.

Since Lemma 4.4 of \cite{Arcs} also holds for a closed subscheme
$S$ of $X \otimes k [t]$ with ${\rm dim}_{k [t]} S < d$,
we see that, for such an $S$, the subset $\cL (S)$ of $\cL (X)$
is a measurable subset of $\cL (X)$ of measure 0. Using Lemma \ref{3.1},
it follows that any $k [t]$-semi-algebraic subset of $\cL (X)$ is
strongly
measurable, with the same measure as in  
\S\kern .15em \ref{pre}.

\medskip

For a measurable subset $A $ of $\cL (X)$ and a function $\alpha : A
\rightarrow
\ZZ \cup \{\infty\}$, we say that $\LL^{-\alpha}$ is {\emph {integrable}}
ot that $\alpha$ is {\emph {exponentially integrable}}
if the fibers of $\alpha$ are measurable and if the motivic integral
$$
\int_{A} \LL^{- \alpha} d\mu
:= 
\sum_{n \in \ZZ} \mu (A \cap \alpha^{-1} (n)) \LL^{- n}
$$
converges in $\widehat \cM$.

\begin{prop}\begin{enumerate}\item[(i)] Let $A_{i}$, $i \in \NN$, be a family of measurable
subsets of $\cL (X)$. Assume the sets $A_{i}$ are mutually disjoint
and that $A := \cup_{i \in \NN}A_{i}$ is measurable.
Then $\sum_{i \in \NN} \mu (A_{i})$
converges in $\widehat \cM$ to $\mu (A)$.
\item[(ii)] If $A$ and $B$ are measurable
subsets of $\cL (X)$ and if $A \subset B$, then
$||\mu (A)|| \leq ||\mu (B)||$.
\end{enumerate}
\end{prop}

\begin{proof}Straightforward exercise, using Lemma \ref{3}.
\end{proof}

\begin{theorem}\label{image}Let $X$ and $Y$ be algebraic varieties
over $k$ of pure dimension $d$, and let $h : \cL (Y) \rightarrow \cL
(X)$
be a $k [t]$-morphism. If $B \subset \cL (Y)$ is measurable,
resp. strongly measurable,
then $h (B) \subset \cL (X)$ is also measurable,
resp. strongly measurable.
\end{theorem}

\begin{proof}We may assume that $Y$ is irreducible. Set
$$
\Delta := 
\cL (Y_{\rm sing})
\cup
h^{-1} (\cL (X_{\rm sing}))
\cup
\Bigl\{y \in \cL (Y) \Bigm | \ord_{t} \cJ_{h} (y) = \infty
\Bigr\}.
$$
We may assume there exists a closed
subscheme $S$ of $Y \otimes k[t]$ with ${\rm dim}_{k[t]}S < d$ such
that $\Delta$ is contained in $\cL (S)$, because otherwise
$h (\cL (Y))$ and $h (B)$ have measure zero.
Since $B$ is measurable and
$\Delta$ is contained in cylindrical subsets $C$ of $\cL (Y)$
with $||\mu (C)||$ arbitrary small, we see that, for every
$\varepsilon > 0$, there exists cylindrical subsets
$B_{i} (\varepsilon)$, $i \in \NN$,
of $\cL (Y)$, such that $B_{0} (\varepsilon) \cap \Delta =\emptyset$,
$B \triangle B_{0} (\varepsilon) \subset \cup_{i \geq 1} B_{i}
(\varepsilon)$,
and $||\mu (B_{i} (\varepsilon))|| < \varepsilon$ for all $i \geq 1$.
Moreover, when $B$ is strongly measurable we can take $B_{0}
(\varepsilon) \subset B$.
Hence $h (B) \triangle h (B_{0} (\varepsilon)) \subset
\cup_{i \geq 1} h(B_{i}(\varepsilon))$. This implies the theorem, since
by Lemma \ref{6.9} below, $h (B_{0} (\varepsilon))$
is cylindrical and, for $i \geq 1$,
$h (B_{i} (\varepsilon))$
is contained in a cylindrical subset
$A_{i} (\varepsilon)$ of $\cL (X)$ with 
$||\mu (A_{i} (\varepsilon))|| \leq \max \, (||\mu (B_{i} (\varepsilon))||,
\varepsilon) \leq \varepsilon$.
\end{proof}

\begin{lem}\label{6.9}Let $X$ and $Y$ be algebraic varieties over $k$,
of pure dimension $d$, and let $h : \cL (Y) \rightarrow \cL (X)$ be a
$k [t]$-morphism.
Let $B$ be a cylindrical subset of $\cL (Y)$. Then the following holds:
\begin{enumerate}
\item[(a)]For every $\varepsilon > 0$,
$h (B)$ is contained in a cylindrical subset $A$ of $\cL (X)$
with $||\mu (A)|| \leq \max \, (||\mu (B)||, \varepsilon)$.
\item[(b)]Assume $B \cap \cL (Y_{\rm sing}) =\emptyset$,
$h (B) \cap \cL (X_{\rm sing}) = \emptyset$, and
$\ord_{t} \cJ_{h} (y)$ is nowhere equal to $\infty$ on $B$.
Then $h (B)$ is cylindrical.
\end{enumerate}
\end{lem}

\begin{proof}(a) First assume that $||\mu (B)|| = 0$. Then, since $B$ is
cylindrical, we have $B \subset \cL (Y_{\rm sing})$ and $h (B)$ is
contained
in some $\cL (S)$, with $S$ a closed subscheme of $X \otimes k[t]$,
with ${\rm dim}_{k[t]}S < d$. This yields
assertion (a) when $||\mu (B)|| = 0$. Now suppose that $||\mu (B)|| \not= 0$.
Take $e$ in $\NN$ large enough to insure that 
$||\mu_{\cL (X)} (\cL (X) \setminus \cL^{(e)} (X))|| \leq ||\mu_{\cL
(Y)} (B)||$. We  may assume that
$h (B)$ is contained in $\cL^{(e)} (X)$. Now we choose
$n \geq e$ 
large enough with respect to $e$ to insure that
$B$ is cylindrical at level $n$ and that
$\cL^{(e)} (Y)$ and 
$\cL^{(e)} (X)$ are cylindrically stable at level $n$.
Set $A := \pi_{n}^{-1} (\pi_{n}(h (B)))$ and note that $A$ is
cylindrical
at level $n$, since $\pi_{n}(h (B))$ is constructible. Moreover, $A$ is
contained in $\cL^{(e)} (X)$, since $h (B)$ is contained in 
$\cL^{(e)} (X)$
and
$n \geq e$. Hence $A$ is cylindrically stable at level $n$.
Thus
\begin{equation}\label{rds}
\mu (A) = [\pi_{n} (h (B))] \, \LL^{- ( n+ 1) d}
\quad
\text{and} 
\quad
||\mu (A) || \leq 2^{- (n  + 1) d + {\rm dim} \pi_{n} (B)}.
\end{equation}
Since
$||\mu (B)|| \not=0$, we have, for $e$ large enough and for $n$ large
enough with respect to $e$, that
$$
{\rm dim} \Bigl(\pi_{n} (B \cap \cL^{(e)} (Y))
\Bigr)
>
{\rm dim} \Bigl(\pi_{n} (\cL (Y) \setminus \cL^{(e)} (Y))
\Bigr) \, ,
$$
and hence $||\mu (B) || = 2^{- (n  + 1) d + {\rm dim} \pi_{n} (B)}$.
Together with (\ref{rds}), this yields assertion (a).

\medskip

(b) Using resolution of singularities, we may assume that $Y$ is smooth.
By
Lemma \ref{3}, there exists $e'$ in $\NN$ such that
$B$ is contained in $h^{-1} (\cL^{(e')} (X))$ and $\ord_{t} \cJ_{h}$
is bounded on $B$. Assertion (b) follows now from the first part of Lemma
\ref{KL}.
\end{proof}

\begin{theorem}[Change of variables formula]\label{CV2}
Let $X$ and $Y$ be algebraic varieties over $k$, of pure
dimension $d$. 
Let $h : \cL (Y) \rightarrow \cL (X )$ be a $k [t]$-morphism and
let 
$A$ and $B$
be 
strongly measurable subsets of $\cL (X)$ and $\cL (Y)$ respectively.
Assume that $h$ induces a bijection between
$B$ and $A$. Then, for any exponentially integrable
function
$\alpha : A \rightarrow
\ZZ \cup \{\infty\}$, the function $B \rightarrow \ZZ \cup \{\infty\} : y \mapsto \alpha (h (y))
+
{\rm ord}_t \cJ_{h} (y)
$ is exponentially integrable and 
$$
\int_A \LL^{- \alpha} d \mu = 
\int_{B}
\LL^{- \alpha \circ h - {\rm ord}_t \cJ_{h} (y)} d \mu.
$$
\end{theorem}

\begin{proof}Reasoning as in the proof of Theorem \ref{image},
we reduce to the case where $B$ is cylindrical and satisfies
$B \cap \Delta = \emptyset$, with $$
\Delta := 
\cL (Y_{\rm sing})
\cup
h^{-1} (\cL (X_{\rm sing}))
\cup
\Bigl\{y \in \cL (Y) \Bigm | \ord_{t} \cJ_{h} (y) = \infty
\Bigr\}.
$$
For this reduction we use the assumption that $B$ is strongly measurable
to insure that the cylinder $B_{0} (\varepsilon)$ in \ref{image}
is contained in $B$, so that the restriction of $h$ to $B_{0}
(\varepsilon)$ is injective.
Next we can reduce to the case where $Y$ is smooth, using resolution of
singularities. Since 
$B \cap \Delta = \emptyset$, it follows from Lemma \ref{3} that there
exists $e'$ in $\NN$ such that $B$ is contained in $h^{-1} (\cL^{(e')}
(X))$
and that 
$\ord_{t} \cJ_{h}$ is bounded on $B$. Thus we may as well assume that
$\ord_{t} \cJ_{h}$ has constant value $e$ on $B$ and the theorem follows
now from Lemma \ref{KL} (b).\end{proof}

\bibliographystyle{amsplain}

\begin{thebibliography}{SGA}



\bibitem{B-N}
S. del Ba{\~ n}o Rollin, V. Navarro Aznar,
\textit{On the motive of a quotient variety},
Collect. Math., \textbf{49} (1998),
203--226.

\bibitem{BD}
V. Batyrev, D. Dais, \textit{Strong McKay correspondence,
string-theoretic Hodge numbers and mirror symmetry}, Topology
\textbf{35}
(1996), 901--929. 

\bibitem{BB}
V. Batyrev, L. Borisov, \textit{Mirror duality and string-theoretic Hodge
numbers}, Inv. Math. \textbf{126} (1996), 183--203.



\bibitem{Ba4}
V. Batyrev,
\textit{Birational Calabi-Yau $n$-folds have equal Betti number},
in {New trends in algebraic geometry (Warwick, 1996)}, {Cambridge Univ. Press},
(1999), 1--11.


\bibitem{Ba2}
V. Batyrev,
\textit{Stringy Hodge numbers of varieties with Gorenstein canonical singularities}, Proc. Taniguchi
Symposium
1997, ``Integrable Systems and Algebraic Geometry, Kobe/Kyoto'',
World Scientific, (1998) 1--32.




\bibitem{Ba}
V. Batyrev,
\textit{Non-archimedian integrals and stringy Euler numbers of log
terminal pairs}, Journal of European Math. Soc. \textbf{1} 
(1999), 5-33.


\bibitem{B-L-R}
S. Bosch, W. L{\"u}tkebohmert, M. Raynaud,
\textit{N{\'e}ron Models},
Ergeb. Math. Grenzgeb. (3) 21, Springer-Verlag, Berlin
1990.




\bibitem{Arcs}
J. Denef, F. Loeser,
\textit{Germs of arcs on singular algebraic varieties
and motivic integration},
Inv. Math.
\textbf{135}
(1999),
201--232.


\bibitem{Motivic}
J. Denef, F. Loeser,
\textit{Motivic Igusa zeta functions},
Journal of Algebraic Geometry,
\textbf{7}
(1998),
505--537.

\bibitem {G-S}
H. Gillet, C. Soul{\'e},
\textit{Descent, motives and $K$-theory},
J. reine angew. Math.
\textbf{478}
(1996),
127--176.

\bibitem {G-N}
F. Guill{\'e}n, V. Navarro Aznar,
\textit{Un crit{\`e}re d'extension d'un foncteur d{\'e}fini sur les
sch{\'e}mas lisses},
preprint (1995), revised (1996).


\bibitem {K}
M. Kontsevich, Lecture at Orsay (December 7, 1995).


\bibitem {Loo}E. Looijenga, {\it Motivic Measures}, in S{\'e}minaire Bourbaki,
expos{\'e} 874, Mars 2000.

\bibitem{Milne}
J. Milne,
\textit{{\'E}tale cohomology},
Princeton University Press (1980),
Princeton, N.J.

\bibitem{Pas}
J. Pas,
\textit{Uniform $p$-adic cell decomposition and local zeta functions},
J. reine angew. Math.
\textbf{399}
(1989),
137--172.

\bibitem{Reid}
M. Reid,
\textit{McKay correspondence},
preprint (1997), available at math. AG/9702016.


\bibitem{Re}M. Reid, {\it La correspondance de McKay}, in S{\'e}minaire
Bourbaki, expos{\'e} 867, Novembre 1999, 20 p., available at math. AG/9911165.

\bibitem{Veys}
W. Veys, 
\textit{The topological zeta function associated to a function on a normal surface
germ}, Topology \textbf{38} (1999), 439--456. 


\end{thebibliography}

\end{document}